\documentclass[11pt]{article}
\usepackage{amsfonts}

\usepackage[colorlinks, linkcolor=blue, citecolor=blue]{hyperref}
\usepackage{epsfig}
\usepackage{epstopdf}
\usepackage{enumerate}
\usepackage{enumitem}
\usepackage{amsmath}
\usepackage{amsfonts}
\usepackage{amssymb}
\usepackage{amsfonts}
\usepackage{amssymb,amsmath,amsthm,epsfig,graphicx,natbib,subfigure}
\usepackage{amssymb,amsmath,graphicx,multirow,booktabs,bm}
\usepackage[margin=1.3in]{geometry}
\usepackage{caption2}
\usepackage{color}
\usepackage{enumerate}
\usepackage{amsfonts}
\usepackage{amssymb,amsmath,amsthm,epsfig,graphicx}
\usepackage{float}
\usepackage{cases}%
\usepackage{tabularx} % set width of a table
\usepackage{multirow}
\usepackage{bigdelim}
\usepackage{color}
\usepackage{graphicx}
\usepackage{soul}
\usepackage{dsfont}  % use the font of mathds
\usepackage{array} % and/or
\usepackage{arydshln} %±í¸ñÐéÏß
\usepackage{longtable} % and/or
\usepackage{colortab} % or
\usepackage{colortbl}
\usepackage{algorithm}

\def\double{\par\baselineskip=22pt}

\def\standard{\oddsidemargin=0in
              \evensidemargin=0in
              \topmargin =-.5in
              \textheight=9.0in
              \textwidth=6.5in}

{\end{list}}
\newenvironment{hangref}{\begin{list}{}{\setlength{\itemsep}{0pt}%
\setlength{\parsep}{0pt}\setlength{\leftmargin}{+\parindent}%
\setlength{\itemindent}{-\parindent}}}{\end{list}}
\def\1{\mathbbm{1}}

\def\E{{\mathbb{E}}}
\def\N{{\mathbb{N}}}
\def\P{{\mathbb{P}}}

\def\R{{\mathbb{R}}}
\def\Var{{\rm Var}}
\def\Cov{{\rm Cov}}

\def\MSE{{\rm MSE}}

\def\ba{\boldsymbol{a}}

\def\bN{{\bf N}}

\def\bt{{\bf t}}

\def\cF{{\cal F}}

\def\cP{{\cal P}}

\def\b1{{\bf 1}}

\standard
\usepackage{authblk}
\usepackage{bbm}
\usepackage{color}
\usepackage{verbatim} %\begin{comment}
\usepackage{threeparttable}

\usepackage{hyperref}

\usepackage{enumitem}
\newtheorem{theorem}{Theorem}
\newtheorem{proposition}{Proposition}%[section]

\newtheorem{definition}{Definition}%[section]
\newtheorem{assumption}{Assumption}
\newtheorem{lemma}{Lemma}
\begin{document}

\normalsize

\title{An Upper Confidence Bound Approach to Estimating the Maximum Mean}

\begin{comment}
\author{Kun Zhang\\
\footnotesize{Institute of Statistics and Big Data, Renmin University of China}\\
\footnotesize{No. 59 Zhongguancun Street, Beijing, China 100872}\\[7pt]
Guangwu Liu\\
\footnotesize{Department of Management Sciences, College of Business,
City University of Hong Kong}\\
\footnotesize{Tat Chee Avenue, Kowloon, Hong Kong}\\[7pt]
Wen Shi,\ Xiaohong Chen\\
\footnotesize{Business School,
Central South University}\\
\footnotesize{Hunan, China 410083}
}
\end{comment}

\author[1]{Kun Zhang\footnote{kunzhang@cityu.edu.hk}}
\author[2]{Guangwu Liu\footnote{msgw.liu@cityu.edu.hk}}
\author[3]{Wen Shi\footnote{shi3wen@163.com}}
\affil[1]{Institute of Statistics and Big Data, Renmin University of China, Beijing, China 100872}
\affil[2]{Department of Management Sciences, College of Business, City University of Hong Kong, Kowloon Tong, Hong Kong}
\affil[3]{Business School, Central South University, Changsha, China 410083}

\date{\today}

\maketitle

\begin{abstract}

Estimating the maximum mean finds a variety of applications in practice. In this paper, we study estimation of the maximum mean using an upper confidence bound (UCB) approach where the sampling budget is adaptively allocated to one of the systems. We study in depth the existing {\it grand average (GA)} estimator, and propose a new {\it largest-size average (LSA)} estimator. Specifically, we establish statistical guarantees, including strong consistency, asymptotic mean squared errors, and central limit theorems (CLTs) for both estimators, which are new to the literature. We show that LSA is preferable over GA, as the bias of the former decays at a rate much faster than that of the latter when sample size increases. By using the CLTs, we further construct asymptotically valid confidence intervals for the maximum mean, and propose a single hypothesis test for a multiple comparison problem with application to clinical trials. Statistical efficiency of the resulting point and interval estimates and the proposed single hypothesis test is demonstrated via numerical examples. 

\end{abstract}

\double

%\section{Format of References}

%This section presents an example which motivates this study, and
%describes the problem to be focused on.  \cite{Artzner99,Artzner99},
%\cite[Section 3.5]{Artzner99} \citeasnoun{AsmBlaJunRo2009},
%\citeaffixed{AsmBlaJunRo2009,Artzner99}{see, e.g.,},
%\citeasnoun[Section 3.5]{AsmGly2007}

%Our scheme is also related to that exact scheme of
%\citeasnoun{AsmBlaJunRo2009}; see also \citeasnoun{Artzner99}.

%Important references: \citeasnoun{KristensenMele11}.

%References: \citeasnoun{FanMarron94}, \citeasnoun{SeifertGasser94},
%\citeasnoun{Wand94} and \citeasnoun{JonesSignorini97}

\section{Introduction}

Estimating the maximum mean of a number of stochastic systems finds a variety of applications in broad areas, ranging from financial risk measurement and evaluation of clinical trials to Markov decision processes and reinforcement learning. In financial risk measurement, coherent risk measures such as expected shortfall play an important role in determining risk capital charge for financial institutions, which is of great practical importance to both risk managers and regulators. In general, a coherent risk measure can be written as the maximum of expected losses (or means) under different probability measures. In practice, its estimation in practice often involves intensive simulation, and developing efficient methods for estimating coherent risk measures has been a topic of interest in simulation community. To improve upon a two-stage procedure in Chen and Dudewicz (1976) for constructing confidence intervals (CIs) for the maximum mean, Lesnevski et al. (2007) proposed a procedure for constructing a two-sided, fixed-width CI for a coherent risk measure, by combining stage-wise sampling and screening based on ranking-and-selection techniques. Estimating the maximum mean also finds applications in evaluation of clinical trials. One of the goals of clinical trials is to estimate the best treatment effect among a number of treatments with unknown mean effects. Estimates of the best mean effect provide useful information into the level of effectiveness of the best treatment.

In both applications in financial risk measurement and clinical trials, the accuracy of the maximum mean estimates is of great concern. For instance, FDA (U.S. Food and Drug Administration) is concerned about the bias and reliability of the estimates of treatment effects when using adaptive designs in clinical trials (FDA 2019). To address such concerns, statistical inference about the maximum mean is highly desirable, providing not only point estimates, but also CIs that offer information about the precision of the estimates. To provide theoretical guarantees for the CIs, central limit theorems (CLTs) are typically needed, which is part of the agenda in our study. In this paper, we aim to develop efficient estimators for the maximum mean, and provide theoretical results on the properties of the estimators, such as their biases, variances, and CLTs. 

CLTs for maximum mean estimators can also be applied to a fundamental hypothesis testing problem in clinical trials, i.e., comparing $K-1$ treatments (with unknown mean effects denoted by $\mu_1,\ldots,\mu_{K-1}$) to a control treatment with known mean effect $\mu_0$, aiming to determine whether any of the $K-1$ treatments outperforms the control; see Dmitrienko and Hsu (2006) and references therein. In the current practice, a commonly used approach is to simultaneously test $K-1$ null hypotheses: $\mu_k\le \mu_0$, $k=1,\ldots,K-1$, where Bonferroni correction is often applied to control the type I error. However, it is well known that this approach may be conservative especially when $K$ is large, and it comes at the cost of increasing the probability of type II errors. With CLTs for estimators of the maximum mean $\mu^* = \max(\mu_0,\mu_1,\ldots,\mu_{K-1})$, we propose an alternative approach by testing only a single hypothesis: $\mu^*=\mu_0$, thus avoiding the use of Bonferroni correction and its subsequent drawbacks. The CLTs provide theoretical guarantees for the validity of the proposed single hypothesis test that is numerically shown to significantly outperform the multiple hypotheses approach. 

Maximum mean estimation has also been studied in the literature on Markov decision processes (MDPs) and reinforcement learning, where the value function for a given state is usually the maximum of multiple expected values taken for different actions. Estimation of the value function is a fundamental sub-task in many solution methods for MDPs and reinforcement learning; see, e.g., Chang et al. (2005), van Hasselt (2010), and Fu (2017). A double (also referred to as cross-validation) estimator was proposed in van Hasselt (2010, 2013) based on two disjoint sets of samples, where one set of samples are used to select the best system (i.e., the system with the maximum mean) and the estimator is given by the average of the selected system using another set of samples. D'Eramo et al. (2016) proposed an estimator by using a weighted average of the sample means of individual systems, where the weight of a system is set as a Gaussian approximation of the probability that this system is the best. The double estimator leads to double Q-learning algorithms, and empirical studies in van Hasselt (2010) and van Hasselt et al. (2016) suggest that these algorithms may outperform conventional Q-learning and deep reinforcement learning algorithms. However, given that maximum mean estimation is only an intermediate module but not the ultimate goal in MDPs and reinforcement learning, how the estimation errors propagate through the state-action pairs and affect the optimal policies is beyond the scope of this paper and requires dedicated future research.

Estimating the maximum mean is closely related to the problem of selecting the best system, an area of active research in simulation community; see, e.g., Kim and Nelson (2006) for an overview. Arguably, the problem of estimating the maximum mean is more difficult than selecting the best system (Lesnevski et al. 2007), because the former requires a careful balance between two goals, including selecting the best system and allocating as many samples as possible to the best system. These two goals are not necessarily complementary, because correctly selecting the best system requires some samples to be allocated to all systems, which limits the achievement of the second goal. 

One of the key issues in estimating the maximum mean is how to sample from individual systems. In this regard, both static sampling (van Hasselt 2010, 2013; D'Eramo et al. 2016) and adaptive sampling (Chang et al. 2005; Kocsis and Szepesv\'ari 2006) have been studied in the literature. The former refers to sampling policies that maintain a pre-specified sampling frequency for each system, while the latter sequentially adjusts sampling frequencies based on historical samples drawn from the systems. Many adaptive sampling methods such as the upper confidence bound (UCB) policy and Thompson sampling share a common feature that higher sampling frequencies are geared towards systems with larger means. When estimating the maximum mean, this feature implies a smaller variance for the resulting estimator because more samples are allocated to the target system, thus offering better computational efficiency. Adaptive sampling/design using policies like UCB offers benefits beyond computational efficiency in the application of clinical trials, where each patient is assigned to one of the $K$ treatments to assess treatment effects. In this application, other than assessing treatment effects, another important goal is to treat patients as effectively as possible during the trials. Therefore, for the well-being of patients, it is desirable to adopt adaptive designs like the UCB policy so that more patients are assigned to more effective treatments. Interested readers are referred to Villar et al. (2015) and Robertson and Wason (2019) for recent work on adaptive designs in clinical trials. 

In this paper, we study maximum mean estimation based on the UCB sampling policy. UCB is one of the most celebrated adaptive sampling methods used in online decision making with learning, which has been extensively used as a tool to balance the trade-off between exploration and exploitation during the course of learning; see, e.g., Agrawal (1995) and Auer et al. (2002) for early work on multi-armed bandit problems, and Hao et al. (2019) and references therein for recent related work. Via the UCB sampling policy, the expected amount of sampling budget allocated to the non-maximum systems is at most of a logarithmic order of the total budget, which implies that the majority of sampling budget is allocated to the system with the maximum mean. %In our study, we consider a generalized UCB policy that allows the exploration rate to take a range of functional forms, and aim to construct efficient estimators of the maximum mean with statistical guarantees.
Within the UCB framework, one of the estimators being considered is the average of all the samples, no matter drawn from the system with the maximum mean or not. We refer to it as the {\it grand average (GA)} estimator. The rationale of the GA estimator stems from the fact that it is dominated by the sample average of the system with the maximum mean, from which the majority of the samples are drawn. The GA estimator is not new in the literature, and has served as a key ingredient for algorithms in MDPs; see, e.g., Chang et al. (2005) and Kocsis and Szepesv\'ari (2006). However, to the best of our knowledge, very little is known about the statistical properties of the GA estimator, except its asymptotic unbiasedness. In this paper, we fill this gap by providing statistical guarantees for the GA estimator, including its strong consistency, mean squared error (MSE), and CLT. On the downside, a well-known drawback of the GA estimator is its bias that results from taking grand average of all systems, including those with smaller means. While this bias gradually dies out as sample size increases, its impact in finite-sample settings cannot be neglected.% especially given that affordable sample sizes may be rather small for large-scale MDPs and reinforcement learning problems.

As an attempt to remedy the above drawback, we propose a new estimator by taking the sample average of the system that has been allocated the largest amount of samples, referred to as the {\it largest-size average (LSA)} estimator. This estimator is intuitively appealing, as the number of allocated samples to the system with the maximum mean dominates those allocated to other systems, ensuring a high probability that the system identified by LSA is exactly the one with the maximum mean. LSA avoids taking average of systems with smaller means, and it is thus expected to have a smaller bias compared to the GA estimator.  %Parallel to the analysis for the GA estimator, we establish strong consistency, CLT, and asymptotic MSE for the LSA estimator. Because it avoids taking average of systems with smaller means, it is expected that the LSA estimator may have smaller bias compared to the GA estimator in the finite-sample setting, which is also illustrated via numerical examples. (now we have theoretical results regarding this)

To summarize, we make the following contributions in this paper.
\begin{itemize}
\item We provide statistical guarantees for the GA estimator, including its strong consistency, MSE, and CLT, which are new to the literature to the best of our knowledge. %Some of the intermediate results in establishing these statistical guarantees may also be important in their own right. In particular, we propose to identify the system with the maximum mean based on ranks of the allocated sampling budgets to individual systems and show its strong consistency which may be useful in developing efficient algorithms for selecting the best system.

\item We develop a new LSA estimator for estimating the maximum mean with statistical guarantees, and show that it is preferable over the GA estimator, due to its bias that converges at least as fast as $\nu_n/n^{3/2}$, a rate much faster than that ($\nu_n/n$) of the GA estimator, where $n$ and $\nu_n$ denote the total sample size and exploration rate, respectively. While polynomial exploration rates are allowed in our setup, both theoretical results and empirical evidence suggest that logarithmic exploration rate is asymptotically optimal when minimizing the MSEs of both the GA and LSA estimators. 

\item Using the established CLTs, we construct asymptotically valid GA and LSA CIs for the maximum mean, and propose a single hypothesis test for a multiple comparison problem with application to clinical trials. Numerical experiments confirm the validity of the CIs, and show that the proposed single test exhibits lower type I error and higher statistical power, compared to a benchmarking test with multiple hypotheses.
\end{itemize}

The rest of the paper is organized as follows. In Section \ref{sec:Backgrounds}, we formulate the problem of estimating the maximum mean. The GA and LSA estimators, and a single hypothesis test for the multiple comparison problem are presented in Section \ref{sec:Estimators}. Theoretical analysis of the estimators is provided in Section \ref{sec:Theory}. A numerical study is presented in Section \ref{sec:NumericalExamples} to demonstrate the efficiency of the estimators and the proposed single hypothesis test, followed by conclusions in Section \ref{sec:Conclusions}.

\section{Problem and Backgrounds}\label{sec:Backgrounds}

Consider $K$ stochastic systems, with their performances denoted by random variables $\{X_k, k=1,\ldots, K\}$. Let $\mu_k=\E(X_k)$. For each $k$, we assume that $X_k$ follows an unknown probability distribution, while independent samples can be drawn. This assumption is usually satisfied in many practical applications. For instance, $X_k$ may represent the stochastic performance of a complex simulation model and follows an unknown probability distribution, while samples of $X_k$ can be generated via simulation. % In learning problems, taking an action may result in a reward $X_k$ that follows an unknown probability distribution, and observations of $X_k$ can be collected during the course of learning.
We make the following assumption on $X_k$'s.
\begin{assumption}\label{ass:subGaussian}
For any $k=1,\ldots,K$, $X_k$ follows a sub-Gaussian distribution, i.e., there exists a sub-Gaussian parameter $\gamma_k$ such that $\E[\exp\left(\lambda(X_k - \mu_k)\right)]\leq \exp\left(\gamma_k^2\lambda^2/2\right)$, for any $\lambda\in\R$.
\end{assumption}

Assumption \ref{ass:subGaussian} is widely used in the literature on machine learning. Many commonly used distributions, such as those with bounded supports and normal distributions, are sub-Gaussian. Assumption \ref{ass:subGaussian} ensures the existence of a Hoeffding bound that serves as the basis for the theoretical analysis in this paper. This bound is presented in Lemma \ref{lem:Hoeffdingsub-Gaussian} in the appendix. 

In this paper, we are interested in estimating the maximum mean of the $K$ systems, defined as
\[\mu^* \triangleq \max_{k=1,\ldots,K}\mu_k.\]
Let $k^*$ denote the index of the system with the largest mean, i.e., $\mu^* = \mu_{k^*}$. Without loss of generality, $k^*$ is assumed to be unique.
\begin{assumption}\label{ass:unique}
The system with the largest mean is unique, i.e., $\mu^*>\mu_k$ for any $k\neq k^*$.
\end{assumption}

Assumption \ref{ass:unique} implies that the mean gap between systems $k^*$ and $k$, given by $\Delta_k\triangleq \mu^*-\mu_k$, is positive for any $k\neq k^*$. 

Given samples of $X_k$, denoted by $\{X_{k,j},j=1,\ldots,n_k\}$, $k=1,\ldots,K$, a straightforward estimator of $\mu^*$, referred to as the {\it maximum estimator} (see, e.g., D'Eramo et al. 2016), is given by
\[\max_{k=1,\ldots,K} \bar X_k,\]
where $\bar X_k[n_k] \triangleq \sum_{j=1}^{n_k}X_{k,j}/n_k$ is the sample average of $X_k$. For notational convenience, where no confusion arises, we write simply $\bar X_k$ in places for $\bar X_k[n_k]$.

It has been well known that the maximum estimator has a positive bias and may overestimate $\mu^*$; see, e.g., Smith and Winkler (2006) and van Hasselt (2010). The positive bias may lead to adverse effects in various applications. For instance, in the estimation of coherent risk measures, it may lead to overestimation of risk, resulting in unduly high capital charges for risky activities (Lesnevski et al. 2007). In reinforcement learning, as errors propagate over all the state-action pairs, the positive bias may negatively affect the speed of learning (van Hasselt 2010). More importantly, when using the maximum estimator, it is not clear how to allocate sampling budget to the systems, i.e., setting the values of $n_k$'s, so as to ensure reasonable accuracy. Therefore, it is of both theoretical and practical interests to develop alternative estimators of $\mu^*$.

\section{An Upper Confidence Bound Approach}\label{sec:Estimators}

Sampling from the $K$ systems is essential for constructing any estimator of $\mu^*$. In developing an efficient estimator, two issues are of major concern. The first issue is on how to efficiently allocate the sampling budget, while the second issue is on methods of constructing an estimator such that statistical guarantees can be established.

To address the first issue, adaptive sampling policies have been studied in the literature. A popular adaptive sampling policy is the UCB policy that was originally proposed for the multi-armed bandit (MAB) problem; see Auer et al. (2002). In traditional MAB problem, in each round, one of the $K$ systems is chosen and a random sample (reward) is drawn from the chosen system, and the decision maker aims to maximize the total rewards collected. At the heart of the MAB problem is to find an adaptive sampling policy that decides which system to sample from in each round so as to balance the tradeoff between exploration and exploitation. Among various sampling policies for MAB, it has been well known that the UCB policy achieves the optimal rate of regret, defined as the expected loss due to the fact that a policy does not always choose the system with the highest expected reward; see Lai and Robbin (1985) and Auer et al. (2002). Specifically, let $T_k(n)$ denote the number of samples drawn from system $k$ during the first $n$ rounds, for $k=1,\dots,K$. The UCB sampling policy is described in Algorithm 1.

\begin{algorithm}[h]
\vspace{12pt} \textbf{Algorithm 1: UCB Sampling Policy}
\vspace{4pt}
 \begin{enumerate}
 \item {\bf Initialization}: During the first $K$ rounds, draw a sample from each system.
 \item {\bf Repeat}: For $n\ge K+1$, draw a sample from the system indexed by
 \begin{equation}
 \kappa_n={\text{arg}\max}_{k\in \{1,\dots,K\}}\left\{\bar X_k[T_k(n-1)]+\sqrt{2\log n\over T_k(n-1)}\right\},\label{eqn:UCB}
 \end{equation}
 where $\bar X_k[T_k(n-1)]$ denotes the sample mean of system $k$ during the first $(n-1)$ rounds.

 Update the sample average of system $\kappa_n$.
 \end{enumerate}
 \end{algorithm}

The UCB policy offers important insights into the adaptive allocation of sampling budget. Essentially, it balances the tradeoff between exploration and exploitation using the UCB defined on the right-hand-side (RHS) of (\ref{eqn:UCB}). On the one hand, systems with higher on-going sample averages have a higher chance to be chosen in the current round, contributing to higher total rewards. On the other hand, systems that are chosen less frequently during the previous rounds, i.e., $T_k(n-1)$ being smaller,  may also have a sufficient chance to be chosen so that such systems will be sufficiently explored. The function $\log n$ in (\ref{eqn:UCB}) can be interpreted as the exploration rate that controls the speed of exploring systems that may not have the largest mean.

In the MAB context, the exploration rate is set to be $\log n$ simply because it leads to the optimal rate of regret. However, it should be pointed out that when the objective is not to minimize the regret, as the case in our problem setting, it is not yet clear whether logarithmic exploration rate is optimal. For the sake of generality, we henceforth allow the exploration rate, denoted by $\nu_n$, to take different forms as a function of $n$, and refer to the resulting UCB policy as a {\it generalized UCB (GUCB) policy}, which is described in Algorithm 2.

\begin{algorithm}[h]
\vspace{12pt} \textbf{Algorithm 2 (GUCB): Generalized UCB Sampling Policy}
\vspace{4pt}
 \begin{enumerate}
 \item {\bf Initialization}: During the first $K$ rounds, draw a sample from each system.
 \item {\bf Repeat}: For $n\ge K+1$, draw a sample from the system indexed by
 \[I_n={\text{arg}\max}_{k\in \{1,\dots,K\}}\left\{\bar X_k[T_k(n-1)]+\sqrt{2\nu_n\over T_k(n-1)}\right\},\]where $\bar X_k[T_k(n-1)]$ denotes the sample mean of system $k$ during the first $(n -1)$ rounds.

 Update the sample average of system $I_n$.
 \end{enumerate}
 \end{algorithm}

\subsection{Grand Average (GA) Estimator}

Before developing estimators for $\mu^*$, we highlight some key properties of the GUCB policy. To convey the main idea, for a while we focus on a special case where the exploration rate $\nu_n = \log n$. In this case, it has been well known that (see Auer et al. 2002), for system $k$ ($k\neq k^*$),
\[\E[T_k(n)]\le C \log n,\]for some constant $C$ that depends on $\Delta_k$. Then, it can be easily seen that $\E[T_{k^*}(n)]$, the expected number of times system $k^*$ is chosen, is of order $n$ during the first $n$ rounds, because the summation of $T_k(n)$'s equals $n$.

In other words, it is expected that among the first $n$ rounds, system $k^*$ is chosen in the majority of rounds. Recall that our objective is to estimate $\mu^*$, the mean of system $k^*$. It is, therefore, reasonable to use the grand average of the samples of all the $n$ rounds as an estimator of $\mu^*$, i.e., the GA estimator. Although it uses the samples that are not drawn from system $k^*$, the validity of the GA estimator can be justified by the fact that the estimation error due to such samples may phase out when $n$ is sufficiently large, as the number of such samples is negligible compared to those drawn from system $k^*$.

Specifically, the GA estimator is given by
\begin{equation}\label{eq: mu_n}
  \widetilde{M}_{n} = \frac{1}{n}\sum_{k=1}^K\sum_{j=1}^{T_k(n)}X_{k,j},
\end{equation}
where $\{X_{k,j}, j=1,\ldots,T_k(n)\}$ denotes the samples drawn from system $k$ during the first $n$ rounds using the GUCB policy in Algorithm 2. 

The GA estimator is not new in the literature, which has served as a key ingredient for algorithms for MDPs and reinforcement learning under the special case that $\nu_n = \log n$. However, research on the statistical properties of the estimator has been underdeveloped. To the best of our knowledge, only the asymptotic unbiasedness of the estimator has been established; see, e.g., Kocsis and Szepesv\'ari (2006), while other statistical properties seem to be missing. One of the contributions of this paper is to fill this gap. In particular, we establish strong consistency, asymptotic MSE, and CLT for the GA estimator, which shall be discussed in detail in Section \ref{sec:Theory}.

By taking grand average of the samples of all the systems including those with means smaller than $\mu^*$, the GA estimator is naturally negatively biased. While this bias gradually dies out as the sample size $n$ increases, its impact cannot be neglected with finite samples, which has been considered as a main drawback of the GA estimator.% When applied to large-scale MDPs and reinforcement learning problems where the affordable sample size $n$ is often small, it has been observed in the literature that the negative impact resulting from the bias of the GA estimator could be substantial; see, e.g., Chang et al. (2005). How to develop estimators with smaller biases has been a topic of both theoretical and practical interest.

\subsection{Largest-Size Average (LSA) Estimator}

As an attempt to remedy the bias drawback of the GA estimator, we propose a new estimator based on a key observation that the number of samples drawn from system $k^*$ usually dominates those from other systems. We therefore define, at the $n$th round,
\begin{equation}\label{eqn:DefInStar}
I_{n}^{*}\triangleq \mathop{\arg\max}_{k=1,...,K} T_k\left(n\right).
\end{equation}
In other words, $I_n^*$ is the index of the system from which the GUCB policy has sampled most frequently so far. It is reasonable to expect that it is very likely to observe $I_n^*=k^*$. Recall that $\mu^*$ is the mean of the performance of system $k^*$. It is, therefore, reasonable to propose a new estimator of $\mu^*$ by taking sample average of the system indexed by $I_n^*$, i.e., the LSA estimator. Specifically, the LSA estimator is given by
\begin{equation}\label{eq: M_n*}
  M_{I_n^*} = \frac{1}{T_{I_n^*}(n)}\sum_{j=1}^{T_{I_n^*}(n)}X_{I_n^*,j},
\end{equation}
where the samples are drawn using the GUCB policy in Algorithm 2. 

The LSA estimator and the maximum estimator share a similar structure, both of them being sample averages of certain identified systems. However, the way of identifying the target system differs significantly. The maximum estimator chooses the system with the largest sample average and uses directly this sample average to estimate $\mu^*$. Instead, the LSA estimator chooses the system from which the GUCB policy has sampled most frequently so far, and uses the sample average of this system to estimate $\mu^*$. The latter seems to be more intuitively appealing in that it is virtually an unbiased sample-mean estimator conditioning on the event that system $k^*$ is correctly identified. Arguably, correctly identifying $k^*$ may be relatively easy as the number of rounds GUCB samples from system $k^*$ is usually overwhelmingly larger than those from other systems. %Compared to the maximum estimator, an advantage of the LSA estimator is its robustness. In the presence of outliers, it is more likely that the maximum estimator turns out to be the sample average of system $k$ with some $k\neq k^*$. By contrast, the LSA is more robust with respect to the outliers, because a few outliers will not affect the identification of the system from which the GUCB policy samples most frequently.

Compared to the GA estimator, the LSA estimator avoids taking averages of those systems with smaller means with a high probability. It is thus expected that the LSA estimator may have smaller bias, which shall be shown theoretically in Section \ref{sec:Theory}.

\subsection{Hypothesis Testing using Maximum Mean Estimators}\label{sec:Test}

In what follows, we demonstrate how the maximum mean estimators may lead to a new approach to a multiple comparison problem that finds important application in clinical trials. 

In a clinical trial, suppose that one is interested in comparing $K-1$ treatments with unknown mean effects $\mu_1,\ldots, \mu_{K-1}$ to a control treatment with known mean effect $\mu_0$, where $K\ge 2$ and a higher mean effect implies a better treatment. In the current practice, a commonly used approach to this comparison problem is to test simultaneously $K-1$ null hypotheses $H_{0,k}$'s against alternatives $H_{1,k}$'s, where
\begin{eqnarray}
H_{0,k}: \mu_k\le \mu_0 \quad v.s. \quad H_{1,k}: \mu_k>\mu_0,\quad k=1,\ldots,K-1.\label{eqn:Temp1Num}
\end{eqnarray}

With the above multiple hypotheses framework, a meta null hypothesis is set as $H_{0,k}$'s being true simultaneously for all $k=1,\ldots,K-1$, which is rejected if any of $H_{0,k}$'s is rejected, implying that at least one of the $K-1$ treatments outperforms the control with certain confidence. Suppose that the overall type I error of this test with multiple hypotheses is set to be $\alpha$. In this case, type I errors associated with individual hypotheses cannot be set as $\alpha$, in order to control the overall type I error. Instead, Bonferroni correction is often applied, which assigns a type I error of $\alpha/(K-1)$ to the $K-1$ individual tests. By doing so, it can be verified that the family-wise error rate (FWER), defined as the probability of rejecting at least one true $H_{0,k}$ is not larger than $\alpha$, provided that these individual tests are independent. To see this, suppose that there are $m_0$ true null hypotheses that are assumed to be $\{H_{0,1},\ldots,H_{0,m_0}\}$ without loss of generality, where $m_0$ takes values in $\{1,\ldots, K-1\}$ and is unknown to the researcher. Then,
\begin{eqnarray*}
\mbox{FWER} = \Pr\left( \cup_{i=1}^{m_0} \{\mbox{Rejecting $H_{0,i}$}\} \right)\le \sum_{i=1}^{m_0}\Pr\left( \mbox{Rejecting $H_{0,i}$} \right)= m_0{\alpha\over K-1}\le \alpha. 
\end{eqnarray*}

While it is easy to implement, Bonferroni correction has been criticized, due to its conservativeness, i.e., a large gap between FWER and $\alpha$ especially when $K$ is large, and the drawback that it usually increases the probability of type II errors. It is thus desirable to develop other procedures for the comparison of treatments that avoid Bonferroni correction. 

We note that with the estimators of the maximum mean, a more appealing approach can be developed for this multiple comparison problem. More specifically, in contrast to the multiple hypotheses formulation, we propose a single hypothesis test, with
\begin{eqnarray}
H_0: \max_{k = 0,\ldots,K-1} \mu_k = \mu_0 \quad v.s. \quad H_1:  \max_{k = 0,\ldots, K-1} \mu_k >\mu_0.\label{eqn:Temp2Num}
\end{eqnarray}

Relabelling $\mu_0$ as $\mu_K$, we note that the GA or LSA estimator, $\widetilde M_n$ or $M_{I_n^*}$, can be used to estimate $ \max_{k = 0,\ldots,K-1} \mu_k \equiv \max_{k = 1,\ldots,K}\mu_k$. Essentially, by comparing $\widetilde M_n$ (or $M_{I_n^*}$) to $\mu_0$, one may be able to decide whether to reject $H_0$ or not. Nevertheless, one needs to know the asymptotic distribution of $\widetilde M_n$ (or $M_{I_n^*}$), in order to formalize the test. We shall show in Section \ref{sec:Theory}, more specifically, Theorems \ref{thm:CLT tilde M} and \ref{thm:CLT Mstar}, that both $\widetilde M_n$ and $M_{I_n^*}$ follow asymptotic normal distributions, thus providing theoretical support for the proposed single hypothesis test. Exact forms of the testing statistics resulting from the GA and LSA estimators are delayed to (\ref{eqn:TestGA}) and (\ref{eqn:TestLSA}) in the subsequent section.

\section{Theoretical Analysis}\label{sec:Theory}

Throughout the paper, we denote $\bar\gamma = \max\{\gamma_1,...,\gamma_K\}$, where $\gamma_k$'s are the sub-Gaussian parameter defined in Assumption \ref{ass:subGaussian}. To facilitate analysis, we first show a result on $T_k(n)$. 

\begin{proposition}\label{prop:mom of Tk}
Suppose that Assumptions \ref{ass:subGaussian} and \ref{ass:unique} hold, and $\nu_n\geq4\bar\gamma^2\log n$. Then, for $k\neq k^*$ and any positive integer $p$,
\begin{eqnarray*}
\E T_k(n)^p\leq \begin{cases}
 \displaystyle \frac{8\nu_n}{\Delta_k^2} + 5 &if\ p=1;\\
 \displaystyle \left\{ \frac{8\nu_n}{\Delta_k^2} + 5 \right\}^2 & if\ p=2;\\
 \displaystyle \left\{\frac{8\nu_n}{\Delta_k^2} + 1 + \left[6\log\left(n+1\right) + O(1)\right]^\frac{1}{3} \right\}^3 & if\ p=3;\\
 \displaystyle \left\{\frac{8\nu_n}{\Delta_k^2} + 1 + \left[\frac{2p}{p-3} (n+1)^{p-3}+O\left(n^{p-4}\right)\right]^\frac{1}{p} \right\}^p &if\ p\geq4,
\end{cases}
\end{eqnarray*}
where $\Delta_k>0$ for $k\neq k^*$ due to Assumption \ref{ass:unique}, and the notation $O(\cdot)$ means that $\limsup_{n\rightarrow\infty} a_n/b_n\le C$ for some constant $C$ if $a_n=O(b_n)$.
\end{proposition}

Proposition \ref{prop:mom of Tk} provides an upper bound for the $p$th moment of $T_k(n)$, the number of samples drawn from system $k$ during the first $n$ rounds, for $k\neq k^*$. This upper bound relies on the exploration rate $\nu_n$, under the assumption that the exploration rate is sufficiently fast, i.e., at least as fast as the logarithmic order. In a special case when $\nu_n=4\bar\gamma^2\log n$ and $p=1$, the result is the same as that in Theorem 1 of Auer et al. (2002).

Essentially, Proposition \ref{prop:mom of Tk} gives an upper bound on the rate at which the $p$th moment of $T_k(n)$ converges to $\infty$ as $n\rightarrow\infty$ for $k\neq k^*$. It implies that the sampling ratios $T_k(n)/n$ may satisfy certain convergence properties. In particular, strong consistency of the sampling ratios are summarized in the following theorem, whose proof is provided in Section \ref{sec:ProofSamplingRatio} of the appendix.

\begin{theorem}\label{thm:Tk a.s.}
Suppose that Assumptions \ref{ass:subGaussian} and \ref{ass:unique} hold, and $\nu_n\in\left[\alpha\log n, n^{1-\delta}\right]$ with $0<\delta<1$ and $\alpha\geq4\bar\gamma^2$. Then, as $n\to\infty$,
\begin{eqnarray*}
\frac{T_k(n)}{n}\overset{a.s.}{\longrightarrow}\begin{cases}
0 &for\ k\neq k^*;\\
1 & for\ k= k^*,
\end{cases}
\end{eqnarray*}
where the notation $\overset{a.s.}{\longrightarrow}$ denotes convergence almost surely (or with probability 1).
\end{theorem}

Theorem \ref{thm:Tk a.s.} implies that in the long run, the number of samples drawn from system $k^*$ dominates the total numbers of samples drawn from other systems, in the almost sure sense. This result suggests that, it is natural to use $I_n^*$, defined in (\ref{eqn:DefInStar}), as an estimator of the unknown $k^*$. Strong consistency of $I_n^*$ can also be established, which is summarized in the following theorem.  Proof of the theorem is provided in Section \ref{sec:ProofInStar} of the appendix.

\begin{theorem}\label{thm: Istar2kstar}
Under the same conditions of Theorem \ref{thm:Tk a.s.}, as $n\to\infty$,
\begin{eqnarray*}
I^*_n \overset{a.s.}{\longrightarrow} k^*.
\end{eqnarray*}
\end{theorem}

\subsection{GA Estimator}\label{sec:GA_Asym}

In what follows, we establish strong consistency, asymptotic MSE and asymptotic normality for the GA estimator $\widetilde M_n$. Let $X_{I_j,j}$ denote the sample drawn at the $j$th round. Note that
\begin{eqnarray*}
\widetilde M_n&=&\frac{1}{n}\sum_{k=1}^K T_k(n) \bar X_k[T_k(n)]=\frac{1}{n}\sum_{j=1}^n X_{I_j,j}=\frac{1}{n}\sum_{j=1}^n \left(X_{I_j,j}-\mu_{I_j}\right)+\frac{1}{n}\sum_{j=1}^n \mu_{I_j}.
\end{eqnarray*}

Define
\begin{eqnarray*}
Z_n\triangleq\sum_{j=1}^n \left(X_{I_j,j}-\mu_{I_j}\right).
\end{eqnarray*}
Then,
\begin{eqnarray}\label{eqn:tilde M decomp}
\widetilde M_n=\frac{Z_n}{n}+\frac{1}{n}\sum_{j=1}^n \mu_{I_j}.
\end{eqnarray}

Let $\cF_n$ be the $\sigma$-field generated by the first $n$ samples for $n\geq1$, and $\cF_0=\{\Omega,\emptyset\}$. Note that $I_n$ is $\cF_{n-1}$-measurable. It follows that for $n\geq1$,
\begin{eqnarray*}
\E\left[\left. X_{I_n,n}\right|\cF_{n-1}\right]=\mu_{I_n},\ {\rm and\ then}\ \E X_{I_n,n}=\E\mu_{I_n}.
\end{eqnarray*}
Therefore, $\E Z_n=0$ and
\begin{eqnarray*}
\E\left[\left. Z_n\right|\cF_{n-1}\right]=Z_{n-1}+\E\left[\left. X_{I_n,n}-\mu_{I_n}\right|\cF_{n-1}\right]=Z_{n-1},
\end{eqnarray*}
i.e., $Z_n$ is a martingale with mean 0.

Note that $X_k-\mu_k$ is sub-Gaussian with parameter $\gamma_k$ for $k=1,...,K$. Because $I_{n+1}$ is $\cF_n$-measurable,
\begin{eqnarray*}
\E\left[\left. \exp\{s(X_{I_{n+1},n+1}-\mu_{I_{n+1}})\}\right|\cF_n\right]\leq \exp\left\{\frac{\bar\gamma^2 s^2}{2}\right\}, \mbox{ for } s\in\R,
\end{eqnarray*}
implies that $X_{I_{n+1},n+1}-\mu_{I_{n+1}}$ is sub-Gaussian with parameter $\bar\gamma$.

Recall that by the definition of $Z_n$,
\begin{eqnarray*}
Z_{n+1} - Z_n = X_{I_{n+1},n+1} - \mu_{I_{n+1}}.
\end{eqnarray*}
Then, we have
\begin{eqnarray*}
\P\left(\left. |Z_{n+1} - Z_n|>a\right|\cF_n\right)\leq 2\exp\left\{-\frac{a^2}{2\bar\gamma^2}\right\},
\end{eqnarray*}
by the Hoeffding bound in Lemma \ref{lem:Hoeffdingsub-Gaussian}.

Further note that $\{Z_{n+1} - Z_n, n=1,2,...\}$ is a martingale difference sequence, the condition of Theorem 2 in Shamir (2011) is thus satisfied. Applying this theorem, for any $\delta>0$, it holds that with probability at least $1-\delta$,
\begin{eqnarray*}
\frac{Z_n}{n} \leq \sqrt{\frac{56\log(1/\delta)}{n/(2\bar\gamma^2)}} = \sqrt{\frac{112\bar\gamma^2\log(1/\delta)}{n}}.
\end{eqnarray*}
Let $\varepsilon=\sqrt{\frac{112\bar\gamma^2\log(1/\delta)}{n}}$. Then, $\delta=\exp\left\{-\frac{n\varepsilon^2}{112\bar\gamma^2}\right\}$ and
\begin{eqnarray*}
\P\left(\frac{Z_n}{n}>\varepsilon\right)\leq \delta = \exp\left\{-\frac{n\varepsilon^2}{112\bar\gamma^2}\right\}.
\end{eqnarray*}
Applying the same argument to the sequence $\{-Z_n, n=1,2,...\}$, symmetrically we have
\begin{eqnarray*}
\P\left(\frac{Z_n}{n}<-\varepsilon\right)\leq \delta = \exp\left\{-\frac{n\varepsilon^2}{112\bar\gamma^2}\right\}.
\end{eqnarray*}
Therefore,
\begin{eqnarray*}
\P\left(\left|\frac{Z_n}{n}\right|>\varepsilon\right)\leq 2\exp\left\{-\frac{n\varepsilon^2}{112\bar\gamma^2}\right\},
\end{eqnarray*}
and moreover,
\begin{eqnarray*}
\sum_{n=1}^\infty\P\left( \left|Z_n/n\right|>\varepsilon\right)< \infty.
\end{eqnarray*}
By Borel-Cantelli lemma, $Z_n/n \overset{a.s.}{\longrightarrow} 0$ as $n\to\infty$.   i.e., the first term on the RHS of (\ref{eqn:tilde M decomp}) converges to 0 almost surely.

Moreover, the second term on the RHS of (\ref{eqn:tilde M decomp}) satisfies
\begin{eqnarray*}
\frac{1}{n}\sum_{j=1}^n \mu_{I_j}=\frac{1}{n}\sum_{k=1}^K T_k(n) \mu_k \overset{a.s.}{\longrightarrow} \mu_{k^*},
\end{eqnarray*}
where the convergence follows from Theorem \ref{thm:Tk a.s.} and the continuous mapping theorem. Therefore, by (\ref{eqn:tilde M decomp}), we establish strong consistency of $\widetilde M_n$ that is summarized as follows.

\begin{theorem}\label{thm:a.s. tilde M}
Under the same conditions of Theorem \ref{thm:Tk a.s.}, $\widetilde M_n$ is a strongly consistent estimator of $\mu^*$, i.e., as $n\to\infty$,
\begin{eqnarray*}
\widetilde M_n\overset{a.s.}{\longrightarrow}\mu^*.
\end{eqnarray*}
\end{theorem}

Theorem \ref{thm:a.s. tilde M} ensures that the GA estimator, $\widetilde M_n$, converges to $\mu^*$ with probability 1, as $n\rightarrow\infty$, which serves as an important theoretical guarantee for the estimator. To further understand its probabilistic behavior, it is desirable to conduct error analysis, especially on its bias and variance. In the following theorem, we establish upper bounds for its bias, variance and MSE. The proof of the theorem is provided in Section \ref{sec:ProofMSE_GA} of the appendix.

\begin{theorem}\label{thm:MSE tilde M}
Suppose that Assumptions \ref{ass:subGaussian} and \ref{ass:unique} hold, and $\nu_n\in\left[\alpha\log n, n^{\frac{1}{2}-\delta}\right]$ with $0<\delta<1$ and $\alpha\geq4\bar\gamma^2$. Let $c$ be a constant satisfying $\frac{8\nu_n}{\Delta_k^2} + 5\leq c\nu_n$ for $k\neq k^*$, and $\bar\mu^2=\max\{\mu_1^2,...,\mu_K^2\}$. Then,
  \begin{eqnarray*}
\left|{\rm Bias}\left(\widetilde M_n\right)\right| \leq\sum_{k\neq k^*}\frac{1}{n}\left( \frac{8\nu_n}{\Delta_k^2}+5\Delta_k\right)=O\left(\frac{\nu_n}{n}\right),\quad \Var\left[\widetilde M_n\right]\le \frac{2K \sigma_{k^*}^2}{n} + \frac{2K^3(\bar\mu^2+1)}{n^2} c^2 \nu_n^2 = O\left(\frac{1}{n}\right),
\end{eqnarray*}
and thus
  \begin{eqnarray*}
{\rm MSE}\left(\widetilde M_n\right)\le {2K \sigma_{k^*}^2\over n} + o\left({1\over n}\right),
\end{eqnarray*}
where the notation $o(\cdot)$ means that $\lim_{n\rightarrow\infty} a_n/b_n=0$ if $a_n=o(b_n)$.
\end{theorem}

Theorem \ref{thm:MSE tilde M} shows that the MSE of the GA estimator converges to $0$ at a rate of $n^{-1}$. In particular, the variance is the dominant term in the MSE, compared to the square of the bias, which is of order $\nu_n^2/n^2$. This result implies that when the exploration rate takes value in the range $\nu_n\in\left[\alpha\log n, n^{1/2-\delta}\right]$ with $0<\delta<1/2$, the bias of the GA estimator is negligible compared to its variance, making it possible to construct CIs by ignoring its bias. Furthermore, as suggested by Theorem \ref{thm:MSE tilde M}, the leading term of the variance does not depend on $\nu_n$. It is thus asymptotically optimal to choose logarithmic $\nu_n$, the slowest possible rate, so as to reduce the asymptotic bias as far as possible. The superiority of this choice of $\nu_n$ is also confirmed by the numerical experiments in Section \ref{sec:NumericalExamples}.

To provide a formal theoretical support for asymptotically valid CIs and the single hypothesis test in (\ref{eqn:Temp2Num}), we establish a CLT for $\widetilde M_n$ in the following theorem, whose proof is provided in Section \ref{sec:ProofCLT_GA} of the appendix.

\begin{theorem}\label{thm:CLT tilde M}
Suppose that Assumptions \ref{ass:subGaussian} and \ref{ass:unique} hold, and $\nu_n\in\left[\alpha\log n, n^{\frac{1}{2}-\delta}\right]$ with $0<\delta<1$ and $\alpha\geq4\bar\gamma^2$. Then,
\begin{eqnarray*}
\sqrt{n}\left(\widetilde M_n-\mu^*\right)\Rightarrow N(0,\sigma_{k^*}^2), \ \ as\ \ n\rightarrow\infty,
\end{eqnarray*}
where $\Rightarrow$ denotes convergence in distribution.
\end{theorem}

Theorem \ref{thm:CLT tilde M} shows that the GA estimator is asymptotically normally distributed with mean $\mu^*$ and variance $\sigma_{k^*}^2/n$. Based on Theorem \ref{thm:CLT tilde M}, an asymptotically valid CI can be constructed for $\mu^*$. To do so, a remaining issue is on how to estimate the unknown $\sigma_{k^*}^2$. In light of the proof of Theorem \ref{thm:CLT tilde M}, we estimate $\sigma_{k^*}^2$ using
\begin{eqnarray}\label{esti var tilde M}
\widetilde\sigma_n^2=\frac{1}{n}\sum_{k=1}^K T_k(n) \hat\sigma_k^2,
\end{eqnarray}
where for $k=1,...,K$,
\begin{eqnarray}\label{esti sigma k}
\hat\sigma_k^2=\frac{1}{T_k(n)}\sum_{j=1}^{T_k(n)}X_{k,j}^2 - \left[\frac{1}{T_k(n)}\sum_{j=1}^{T_k(n)}X_{k,j}\right]^2.
\end{eqnarray}

It can be shown that $\widetilde\sigma_n^2$ converges to $\sigma_{k^*}^2$ in probability, as $n\rightarrow\infty$. This result is summarized in the following proposition, whose proof is provided in Section \ref{sec:ProofSigmaConsistency} of the appendix.
\begin{proposition}\label{prop:CI tilde M}
Under the same conditions of Theorem \ref{thm:CLT tilde M}, as $n\rightarrow\infty$,
\begin{eqnarray*}
\widetilde\sigma_n^2\to \sigma_{k^*}^2,\ \ {\rm in\ probability}.
\end{eqnarray*}
\end{proposition}

From Theorem \ref{thm:CLT tilde M} and Proposition \ref{prop:CI tilde M}, it follows that
\begin{eqnarray*}
\sqrt{n} \widetilde\sigma_n^{-1}\left(\widetilde M_n-\mu^*\right)\Rightarrow N(0,1).
\end{eqnarray*}
Then, an asymptotically valid 100$(1-\beta)\%$ CI of $\mu^*$ is given by
\begin{eqnarray}\label{CI tilde}
\left(\widetilde M_n-z_{1-\beta/2}\widetilde\sigma_n/\sqrt{n}, \quad \widetilde M_n+z_{1-\beta/2}\widetilde\sigma_n/\sqrt{n}\right),
\end{eqnarray}
where $z_{1-\beta/2}$ is the $1-\beta/2$ quantile of the standard normal distribution.

The CLT in Theorem \ref{thm:CLT tilde M} and the variance estimator $\widetilde \sigma_n^2$ in Proposition \ref{prop:CI tilde M} also serve as the theoretical foundation for the proposed single hypothesis test in Section \ref{sec:Test}. More specifically, the following testing statistics can be constructed:
\begin{equation}\label{eqn:TestGA}
\widetilde U_n \triangleq {\sqrt{n}\left(\widetilde M_n - \mu_0\right) \over \widetilde \sigma_n}.
\end{equation}
Note that under either the null or alternative hypothesis, $\max_{k=0,\ldots,K-1} \mu_k \ge \mu_0$. Therefore, a one-side test can be conducted using the GA-based statistics $\widetilde U_n$. In particular, the null hypothesis is rejected when $\widetilde U_n > z_{1-\alpha}$, suggesting that at least one of the considered treatments outperform the control treatment, where $\alpha$ denotes the overall type I error.

\subsection{LSA Estimator}\label{sec:LSA_Asym}

Strong consistency of the LSA estimator is stated in the following theorem, whose proof is provided in Section \ref{sec:ProofConsistencyLSA} of the appendix.
\begin{theorem}\label{TIstar a.s.}
Suppose that Assumptions \ref{ass:subGaussian} and \ref{ass:unique} hold, and $\nu_n\in\left[\alpha\log n, n^{1-\delta}\right]$ with $0<\delta<1$ and $\alpha\geq4\bar\gamma^2$. Then, as $n\to\infty$,
\begin{eqnarray*}
M_{I^*_n}\overset{a.s.}{\longrightarrow}\mu^*.
\end{eqnarray*}
\end{theorem}

Theorem \ref{TIstar a.s.} shows that the LSA estimator converges to $\mu^*$ as $n\rightarrow\infty$. To better understand its asymptotic error, we establish upper bounds for its asymptotic bias, variance and MSE, which are summarized in the following theorem with proof provided in Section \ref{sec:ProofMSE_LSA} of the appendix.

\begin{theorem}\label{thm:MSE Mstar}
Under the same conditions of Theorem \ref{TIstar a.s.}, ${\rm Bias}(M_{I_n^*}) = \E\left[M_{I_n^*} - \mu_{k^*}\right]<0$ for sufficiently large $n$, and there exist constants $c_1,c_2$ and $c_3$, such that
\begin{align*}
\left|{\rm Bias}\left( M_{I_n^*} \right)\right| \le {K^2\sigma_{k^*}^2 c_2\nu_n + c_3\over \sqrt{n}\left( n - Kc_1\nu_n \right)} +o\left( {\nu_n\over n^{3/2}} \right) = O\left( {\nu_n\over n^{3/2}} \right),
\end{align*}
\begin{align*}
\Var\left[M_{I_n^*}\right] \leq {\sigma_{k^*}^2 K^2\over n} + o\left( {1\over n}\right),
\end{align*}
and thus
\begin{align*}
\MSE\left[M_{I_n^*}\right] \leq{\sigma_{k^*}^2 K^2\over n} + O\left(\frac{\nu_n^2}{n^3}\right) = O\left(\frac{1}{n}\right).
\end{align*}
\end{theorem}

Theorem \ref{thm:MSE Mstar} shows that the bias of LSA estimator decays at a rate at least as fast as $\nu_n/n^{3/2}$, which is much faster than that of GA estimator (see Theorem \ref{thm:MSE tilde M}). Interestingly, the rate of convergence of its bias can be faster than $n^{-1}$ if $\nu_n=o(n^{1/2})$. Furthermore, as suggested by Theorem  \ref{thm:MSE Mstar}, the leading term of the variance does not depend on $\nu_n$. It is thus asymptotically optimal to choose logarithmic $\nu_n$, the slowest possible rate, so as to reduce the asymptotic bias as far as possible. The superiority of this choice of $\nu_n$ is also confirmed by the numerical experiments in Section \ref{sec:NumericalExamples}.

It should be remarked that in Theorems \ref{thm:MSE Mstar} and \ref{thm:MSE tilde M}, the dominant coefficients of the variance bounds for the LSA and GA estimators, i.e., $\sigma_{k^*}^2K^2$ and $2K\sigma_{k^*}^2$, are not tight and can be further reduced with more elaborate analysis. These coefficients do not imply that the GA estimator has a smaller variance, and in fact, numerical results in Section \ref{sec:NumericalExamples} suggest the opposite. 

Similar to the case for GA, we establish a CLT for the LSA estimator, which is stated in the following theorem. Proof of the theorem is provided in Section \ref{sec:ProofCLT_LSA} of the appendix.

\begin{theorem}\label{thm:CLT Mstar}
Under the same conditions of Theorem \ref{TIstar a.s.},
\begin{eqnarray*}
\sqrt{n}\left(M_{I^*_n}-\mu^*\right) \Rightarrow N(0,\sigma_{k^*}^2),\ \ as\ \ n\to\infty.
\end{eqnarray*}
\end{theorem}

Theorem \ref{thm:CLT Mstar} shows the asymptotic normality of the LSA estimator, whose rate of convergence is $n^{-1/2}$. From its proof, it can also be seen that if replacing $\sqrt{n}$ by $\sqrt{T_{I^*_n}(n)}$, the CLT still holds.

While both CLTs for the GA and LSA estimators (in Theorems \ref{thm:CLT tilde M} and \ref{thm:CLT Mstar}, respectively) have the same rate of convergence, it should be emphasized that they require different conditions on the exploration rate. Specifically, CLT holds for the LSA estimator so long as the exploration rate $\nu_n$ is within the range of $\left[\alpha\log n, n^{1-\delta}\right]$, while CLT of the GA estimator can only be guaranteed under a more restricted condition that $\nu_n\in \left[\alpha\log n, n^{{1/2}-\delta}\right]$.

In what follows, we construct an asymptotically valid CI for $\mu^*$ based on Theorem \ref{thm:CLT Mstar}. To do this, we propose to estimate the unknown variance $\sigma_{k^*}^2$ by
\begin{eqnarray*}
\bar\sigma_{n}^2=\frac{1}{T_{I^*_n}(n)}\sum_{j=1}^{T_{I^*_n}(n)}X_{I^*_n,j}^2-\left[\frac{1}{T_{I^*_n}(n)}\sum_{j=1}^{T_{I^*_n}(n)}X_{I^*_n,j}\right]^2.
\end{eqnarray*}

Strong consistency of the estimator $\bar\sigma_n^2$ is established in the following proposition, whose proof is provided in Section \ref{sec:ProofConsistencyVar_LSA} of the appendix.
\begin{proposition}\label{prop:CI Mstar}
Under the same conditions of Theorem \ref{TIstar a.s.},
\begin{eqnarray}\label{bar sigma}
\bar\sigma_n^2\overset{a.s.}{\longrightarrow} \sigma_{k^*}^2,\quad {\mbox as} \ n\rightarrow\infty.
\end{eqnarray}
\end{proposition}

Similar to (\ref{CI tilde}), an asymptotically valid 100$(1-\beta)\%$ CI of $\mu^*$ is given by
\begin{eqnarray}\label{CI star}
\left(M_{I^*_n}-z_{1-\beta/2}\bar\sigma_n/\sqrt{n}, \quad M_{I^*_n}+z_{1-\beta/2}\bar\sigma_n/\sqrt{n}\right).
\end{eqnarray}

Similar to the discussion in Section \ref{sec:GA_Asym}, for the proposed single hypothesis test in Section \ref{sec:Test}, we construct the following LSA-based testing statistics:
\begin{equation}\label{eqn:TestLSA}
U_n^* \triangleq {\sqrt{n}\left(M_{I_n^*} - \mu_0\right) \over \bar \sigma_n},
\end{equation}
and reject the null hypothesis when $U_n^* > z_{1-\alpha}$, where $\alpha$ denotes the overall type I error.

\section{Numerical Study}\label{sec:NumericalExamples}

We consider three examples. The first example aims to illustrate the asymptotic behaviors of the GA and LSA estimators, including their errors and rates of convergence. While it has been argued in Sections \ref{sec:GA_Asym} and \ref{sec:LSA_Asym} that logarithmic exploration rate is asymptotically optimal when minimizing MSEs of the GA and LSA estimators, we numerically confirm it using the first set of experiments. For the rest of the experiments in all examples, we choose logarithmic exploration rate given its superiority on both theoretical and empirical grounds. 

In the second example, we consider an application of the GA and LSA estimator to the estimation of coherent risk measures, and compare our constructed GA and LSA CIs to the fixed-width CIs produced by the procedure of Lesnevski et al. (2007).

In the third example, we consider a clinical trial using the GA- and LSA-based single hypothesis tests proposed in (\ref{eqn:Temp2Num}), aiming to demonstrate their relative merits compared to an existing test with multiple hypotheses in the form of (\ref{eqn:Temp1Num}).

\subsection{Implementation Issues}

In the numerical experiments, we observed that the bias of the GA estimator with finite samples could be quite large, especially when the variances of the systems are large, or the exploration rate $\nu_n$ is large. This phenomenon is intuitively explained by that there may be too many samples drawn from systems other than $k^*$, leading to a significant underestimation of the maximum mean. However, this adverse effect dies out gradually in the long run. As a means to alleviate this adverse effect in finite-sample setting, we propose to choose a certain amount of samples as {\it a warm-up period} and discard the warm-up samples in the estimation. Numerical experiments suggest that a small amount of warm-up samples, e.g. $10\%$ of the total sampling budget, may greatly improve the estimation quality of the GA estimator. Throughout our experiments, we set 10\% of the budget as the warm-up period for the GA estimator.

By contrast, the LSA estimator is insensitive to the warm-up samples. Whether to keep or to discard these samples does not have much impact on the quality of the LSA estimator. During the implementation, we discard the warm-up period samples for the LSA estimator as well, to be parallel to the setting for the GA estimator.

From an implementation perspective, it has also been documented that slightly tuning the UCB policy by taking into consideration of the variances of the systems may substantially improve the performances of algorithms for multi-armed bandit problems; see, e.g., Auer et al. (2002). We observe the same phenomenon when estimating the maximum mean. Following the same idea of Auer et al. (2002), we propose to replace the UCB, i.e., $\sqrt{2\nu_n/ T_k(n-1)}$, by
\[\sqrt{{2\nu_n\over T_k(n-1)} V_{k,n}}\]for $k=1,\ldots,K$ during the implementation, where
\[V_{k,n}=\frac{1}{T_k(n-1)}\sum_{i=1}^{T_k(n-1)}X_{k,j}^2 - \bar X_k^2[T_k(n-1)] + \sqrt{2\nu_n\over T_k(n-1)}\]is an estimate of the UCB of the variance of system $k$.

\subsection{Maximum of $K$ Normal Means}\label{sec:Example1}

Consider $K$ systems, with samples from normal distributions with unknown means $(\mu_1, \mu_2, \dots, \mu_K)$. Parameter configuration is similar to that of Fan et al. (2016), who focused on selecting the system with the largest mean rather than estimating the maximum mean. More specifically, the means $(\mu_1, \mu_2, \dots, \mu_K)$ of systems are set to be  $\mu_k = 1.5 + 0.5k$, $k=1,\ldots,K$, and the standard deviations are set to be $\sigma_k = \mu_k$ for all systems.

We examine the performances of different estimators when estimating $\mu^*$, by considering their biases, standard deviations, MSEs, and RRMSEs, where RRMSE is defined as the percentage of the root MSE to the true maximum mean. Results reported are based on $10^4$ independent replications unless otherwise stated.

\begin{figure}[hptb]
	\begin{center}
         \hspace*{-1.78cm}
        \includegraphics[scale=0.27]{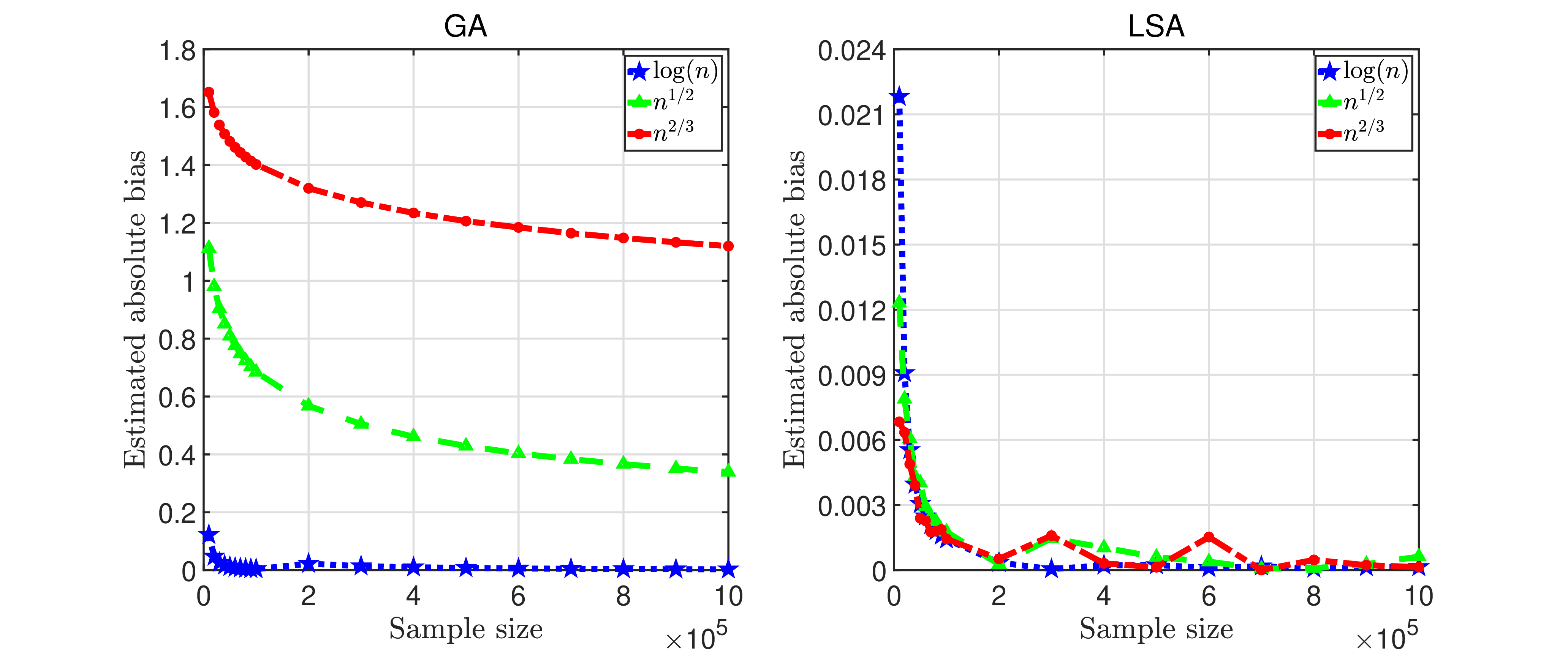}
		\caption{(Color online) Convergence of estimated absolute biases.}
		\label{fig:BiasCompare}
	\end{center}
\end{figure}

\begin{figure}[hptb]
	\begin{center}
         \hspace*{-1.78cm}
        \includegraphics[scale=0.288]{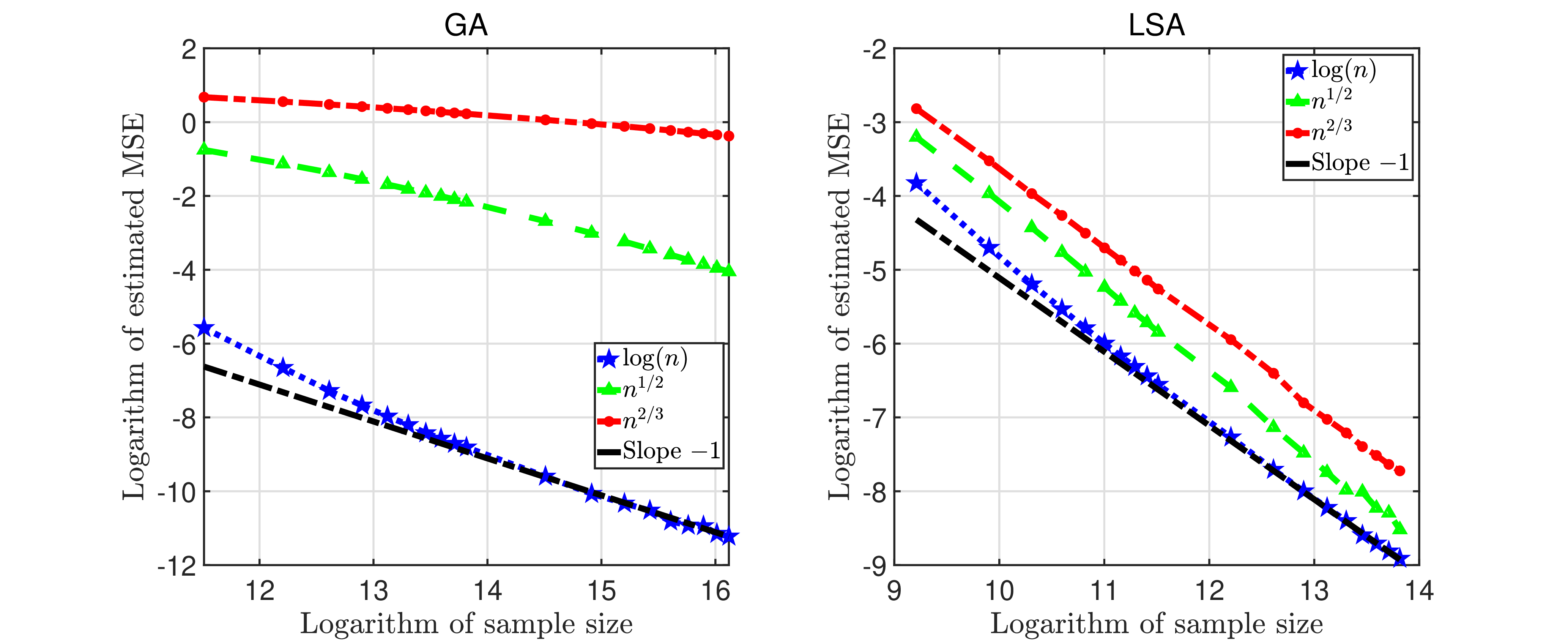}
		\caption{(Color online) Convergence rate of estimated MSEs.}
		\label{fig:MSERate_mu}
	\end{center}
\end{figure}

\begin{figure}[hptb]
	\begin{center}
         \hspace*{-1.78cm}
        \includegraphics[scale=0.2878]{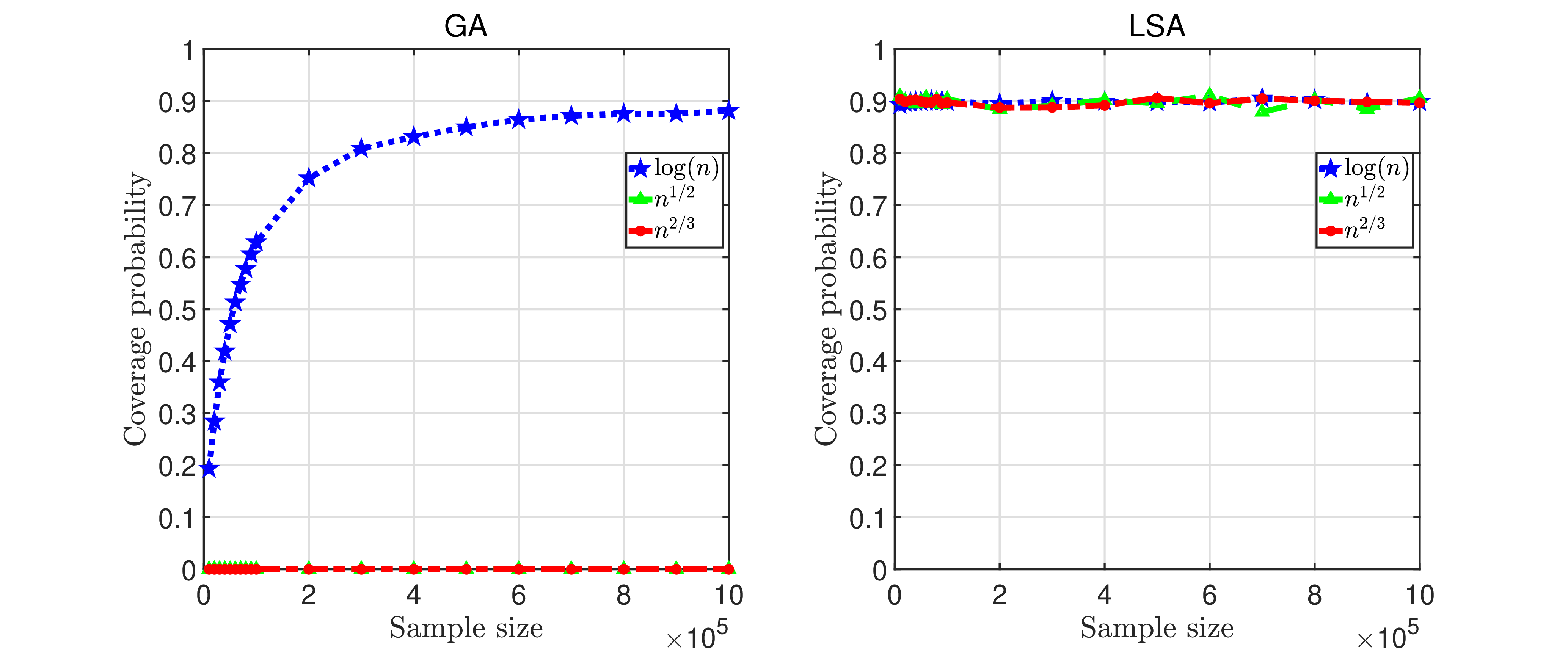}
		\caption{(Color online) Observed coverage probabilities of the 90\% confidence intervals with respect to different sample sizes.}
		\label{fig:CovProbs_mu}
	\end{center}
\end{figure}

In the first set of experiments, let $K=20$ with the total sampling budget $n$ ranging between $10^4$ and $10^7$. We compare biases, MSEs, and coverage probabilities for GA and LSA estimators with respect to (w.r.t.) $\nu_n = \log n$, $n^{1/2}$, and $n^{2/3}$. The results are summarized in Figures \ref{fig:BiasCompare}, \ref{fig:MSERate_mu}, and \ref{fig:CovProbs_mu}, illustrating the rates of convergence of absolute biases and MSEs, and coverage probabilities of the 90\% CIs, for the GA (left panel) and LSA (right panel) estimators.

\begin{itemize}

  \item Figure \ref{fig:BiasCompare} shows that for different choices of $\nu_n$, the bias of the LSA estimator is substantially smaller than that of the GA estimator, consistent with Theorems \ref{thm:MSE tilde M} and \ref{thm:MSE Mstar} that the former converges to $0$ at a rate of at least $\nu_n/n^{3/2}$, much faster than the latter (with a rate of $\nu_n/n$). Furthermore, it is observed that logarithmic $\nu_n$ leads to the smallest biases among the three choices, for both the GA and LSA estimators. 

  \item Figure \ref{fig:MSERate_mu} shows that the MSE of the LSA estimator converges approximately to $0$ at a rate of $n^{-1}$ for all choice of $\nu_n$. For the GA estimator, the same rate is observed for logarithmic $\nu_n$, while the rate tends to be slower when $\nu_n=n^{1/2}$ and $\nu_n=n^{2/3}$. This coincides with the requirement in Theorems \ref{thm:MSE tilde M} that $\nu_n$ must be slower than $n^{1/2}$ for the GA estimator.
  
  \item Figure \ref{fig:CovProbs_mu} shows the LSA CIs outperform the GA ones substantially in terms of coverage probabilities, especially when sample sizes are relatively small, due to the significantly smaller bias of the LSA estimator. When $\nu_n=\log n$, both the GA and LSA CIs achieve the nominal coverage probability with sufficiently large sample sizes. However, the GA CIs fail to do so for $\nu_n=n^{1/2}$ and $\nu_n=n^{2/3}$ that violate the requirement in Theorem \ref{thm:CLT tilde M}. 
\end{itemize}

Given the superiority of logarithmic exploration rate, we set $\nu_n=\log n$ in the rest of the experiments in Section \ref{sec:NumericalExamples}. 

In the second set of experiments, we fix $K=20$, and compare the GA and LSA estimators to four competing estimators given as follows. Consider two maximum average estimators that take the maximum of the sample means of the $K$ systems. In particular, the first is a maximum average estimator based on the GUCB adaptive sampling policy, referred to as adaptive maximum average (AMA) estimator, while the other is the maximum average estimator based on a simple static sampling policy that randomly allocates each sample to one of the systems with equal probabilities, referred to as simple maximum average (SMA) estimator. Moreover, we consider two heuristic estimators in the form of the sample average of certain selected system resulting from existing best arm identification methods. In particular, we consider two Bayesian best arm identification methods, including the top-two probability sampling (TTPS) method of Russo (2020), and the top-two expected improvement (TTEI) of Qin et al. (2017). For these two estimators, we assume known variances of the systems as required by the TTPS and TTEI methods, and set the tuning parameter $\beta$ to be the default value of $0.5$, following Russo (2020) and Qin et al. (2017).

It should be pointed out that, in spite of being reasonable estimators, asymptotic properties of the AMA, TTPS and TTEI estimators, such as their biases, MSEs and central limit theorems, are not known in the literature. The comparison of these estimators in our experiments serves merely to numerically examine their merits, while development of theoretical underpinnings is a topic that deserves future investigation.

Relative biases and standard deviations as percentages of the true maximum mean, and RRMSEs are summarized in Table \ref{tab: Compare}\footnote{For cases of $n=10^4$ and $10^5$, reported errors are estimated based on $10^5$ independent replications.}. From Table \ref{tab: Compare}, it can
be seen that all the considered estimators converge as sample size increases. Among them, LSA and AMA have almost the same biases and variances, and perform the best with the smallest MSEs. Their biases are also smaller or comparable to other estimators. GA and SMA suffer from a disadvantage of significantly larger biases compared to other estimators. TTPS and TTEI have biases comparable to those of LSA and AMA, but with larger variances. It is also interesting to observe that while SMA is positively biased, AMA is negatively biased, suggesting that the dependence structure among the samples drawn from the UCB policy may lead to a negative bias.

\begin{table}[hptb]%\footnotesize
\renewcommand\arraystretch{1.2}
\centering
\caption{Relative bias (\%), relative standard deviation (\%) and RRMSE (\%) of different estimators}
\vspace{0.1in}
\label{tab: Compare}
\begin{tabular}{crrrrrrrrrr}
\hline
 %{\rule[-1mm]{0mm}{17mm}}
          & \multicolumn{3}{c}{LSA}&&\multicolumn{3}{c}{GA}\\
             \cline{2-4}\cline{6-8}
$n$                 &$10^4$   &$10^5$        &$10^6$       &&$10^4$   &$10^5$    &$10^6$\\
\hline
 Relative bias      &$-$0.12  &$-$2.7e-03 &$-$1.3e-03 &&$-$2.38 &$-$0.40  &$-$0.03 \\
 Relative stdev     &1.28     &0.33        &0.10        &&1.31    &0.36     &0.10  \\
 RRMSE              &1.28     &0.33        &0.10        &&3.03    &0.54     &0.11  \\
\hline
          & \multicolumn{3}{c}{AMA}&&\multicolumn{3}{c}{SMA}\\
             \cline{2-4}\cline{6-8}
$n$                 &$10^4$   &$10^5$        &$10^6$       &&$10^4$   &$10^5$    &$10^6$\\
\hline
 Relative bias      &$-0.12$ &$-$2.7e-03  &$-$1.3e-03 &&0.98    &0.01     &3.0e-05 \\
 Relative stdev     &1.27    &0.33         &0.10        &&3.76    &1.40     &0.45   \\
 RRMSE              &1.27    &0.33         &0.10        &&3.89    &1.40     &0.45   \\
\hline
          & \multicolumn{3}{c}{TTPS}&&\multicolumn{3}{c}{TTEI}\\
             \cline{2-4}\cline{6-8}
$n$                 &$10^4$   &$10^5$        &$10^6$       &&$10^4$   &$10^5$    &$10^6$\\
\hline
 Relative bias   &$-$9.0e-03 &$-$7.3e-03  &$-$6.1e-03 && $-0.14$ &$-$2.7e-03 &$-$7.5e-04\\
 Relative stdev     &1.51    &0.49         &0.14        &&1.77    &0.45     &0.14   \\
 RRMSE              &1.51    &0.50         &0.14        &&1.77    &0.45     &0.14   \\
\hline
\end{tabular}
\end{table}

\subsection{Simulated CIs for Coherent Risk Measures}\label{sub: Large-Scale}

Consider a risk measurement application in which we are interested in estimating a coherent risk measure that is the maximum of the expected discounted portfolio loss taken over $K$ probability distributions. Lesnevski et al. (2007) proposed a simulation-based procedure to construct a fixed-width CI for the risk measure, and used an example with $K=256$ to demonstrate the performance of the constructed CIs. For the sake of completeness, we briefly describe this example as follows.

A coherent risk measure can be written as
\begin{eqnarray*}
\sup_{\P\in\cP} \E_\P [-Y/r],
\end{eqnarray*}
where $Y$ is the value of a portfolio at a future time horizon, $1/r$ is a stochastic discount factor which represents the time value of money, and $\cP$ is a set of probability measures simplified by $\cP=\{\P_1,...,\P_K\}$. Denote $X=-Y/r$ and $\mu_k=\E_{\P_k}[X]$. Let $X_i$ be a random variable whose distribution under a probability measure $\P$ is the same as that of $X$ under $\P_k$, i.e., $\P(X_k\leq x) = \P_k(X\leq x)$. Then the coherent risk measure is exactly the maximum mean:
\begin{eqnarray}\label{coherent}
\max_{k=1,...,K} \E X_k.
\end{eqnarray}

Consider estimating the coherent risk measure of a portfolio consisting of European-style call and put options written on three assets with price $S_j(t)$ at time $t$ for $j=1,2,3$. Let $S_0(t)$ denote a market index. All options in the portfolio have a common terminal time $T$. For each of $j=0,1,2,3$, $S_j(t)$ follows geometric Brownian motion with drift $d_j$ and volatility $\zeta_j$, and therefore $\log S_j(t) = \log S_j(0) + (d_j-\zeta_j^2/2)t + \zeta_j W_j \sqrt{t}$, where $W_j$ is a standard normal random variable. The assets are dependent in the sense that $W_0=Z_0$, and $W_j=\lambda_j Z_0+\sqrt{1-\lambda_j^2} Z_j$ for $j=1,2,3$, where under probability measure $\P$, $Z_0$, $Z_1$, $Z_2$, and $Z_3$ are independent standard normal random variables. In this model, $Z_0$ corresponds to the market factor common to all assets, while $Z_1$, $Z_2$, and $Z_3$ are idiosyncratic factors corresponding to each individual asset.

Denote by $\Phi$ the standard normal distribution function. Each $Z_j$ for $j=0,1,2,3$ is sampled restricted on four cases with equal probabilities: sampled conditional on exceeding $\Phi^{-1}(1-1/\sqrt{20})$, sampled conditional on being below $\Phi^{-1}(1/\sqrt{20})$, sampled conditional on being between these two numbers, and sampled with no restriction. Therefore, $K=4^4=256$ in (\ref{coherent}), and each $X_k$ denotes the value of the portfolio corresponding to a combination of these four cases. The parameters are set as follows: $T$ is one week, $\zeta_1=39.8\%$, $\zeta_2=19.3\%$, $\zeta_3=27.0\%$, $\lambda_1=0.617$, $\lambda_2=0.368$, $\lambda_3=0.785$, and $d_j=0$ and $S_j(0)=100$ for $j=1,2,3$. The amounts of options in the portfolio could be referred to Table 1 in Lesnevski et al. (2007).

We estimate the coherent risk measure using the GA and LSA estimators. In the experiments, we examine the performance of the estimators, and compare the CIs to those of Lesnevski et al. (2007), which are constructed via the so-called ``two-stage algorithm with screening" in the literature. When implementing this algorithm, we use exactly the same algorithm parameters as in Lesnevski et al. (2007). %During the implementation of the GA and LSA estimators, we observed that their performances are robust w.r.t. different exploration rates, and the reported results are based on an exploration rate $\nu_n=\log n$.

The simulation-based procedure of Lesnevski et al. (2007) sets a certain amount of initial samples to enable the use of ranking-and-selection techniques. In our implementation, we set the same amount of warm-up  samples. Moreover, a fixed width of the CI has to be set at the beginning of their procedure, the total sample size required to stop the process may thus vary across different replications. For the sake of fairness in the comparison, we first run the procedure of Lesnevski et al. (2007), and the same sample size is then used to construct the GA and LSA CIs. We then compare the width of our CIs to that of Lesnevski et al. (2007).

Relative errors of the GA and LSA estimators and coverage probabilities of their $95\%$ CIs w.r.t. different sample sizes are presented in Table \ref{tab: Nelson paper}. From the table, it can be seen that both the GA and LSA estimators exhibit high accuracy. In particular, when the sample size is larger than $10^5$, RRMSEs of both estimators are less than $1\%$, which are contributed mainly by the variances while their biases are negligible. However, when the sample size is relatively small, e.g., $n=10^4$, the bias of the LSA estimator is substantially lower than that of the GA estimator, which also explains the low coverage probabilities of the GA CIs in this case. This result suggests that the LSA estimator is preferable over the GA estimator, especially for small or medium sample sizes.

\begin{table}[]%\footnotesize
\renewcommand\arraystretch{1.2}
\centering
\caption{Relative bias (\%), standard deviation (\%) and RRMSE (\%), and coverage probabilities (\%)}
\vspace{0.1in}
\label{tab: Nelson paper}
\begin{tabular}{crrrrrrrrr}
%{p{2.5cm}<{\centering}p{0.8cm}<{\centering}p{0.8cm}<{\centering}p{0.8cm}<{\centering}p{0.8cm}<{\centering}p{0.001cm}p{0.7cm}<{\centering}p{0.7cm}<{\centering}p{0.7cm}<{\centering}p{0.8cm}<{\centering}}
\hline
 %{\rule[-1mm]{0mm}{17mm}}
          & \multicolumn{4}{c}{LSA}&&\multicolumn{4}{c}{GA}\\
             \cline{2-5}\cline{7-10}
$n$                 &$10^4$   &$10^5$    &$10^6$  &$10^7$ &&$10^4$  &$10^5$  &$10^6$  &$10^7$\\
\hline
 Bias               &$-0.07$  &$-0.01$   &$-0.00$ &$0.00$ &&$-10.64$&$-0.48$ &$-0.03$ &$-0.00$\\
 Stdev              &1.21	  &0.29	     &0.09	  &0.03   &&1.54    &0.31    &0.09    &0.03\\
 RRMSE              &1.21	  &0.29	     &0.09	  &0.03   &&10.75   &0.57    &0.10    &0.03\\
 Cover. prob.       &93.1     &95.7      &95.3    &95.8   &&0       &60.6    &93.0    &95.6\\
\hline
\end{tabular}
\end{table}

\begin{figure}[ptb]
	\begin{center}
         \hspace*{-1.17cm}
        \includegraphics[scale=0.55]{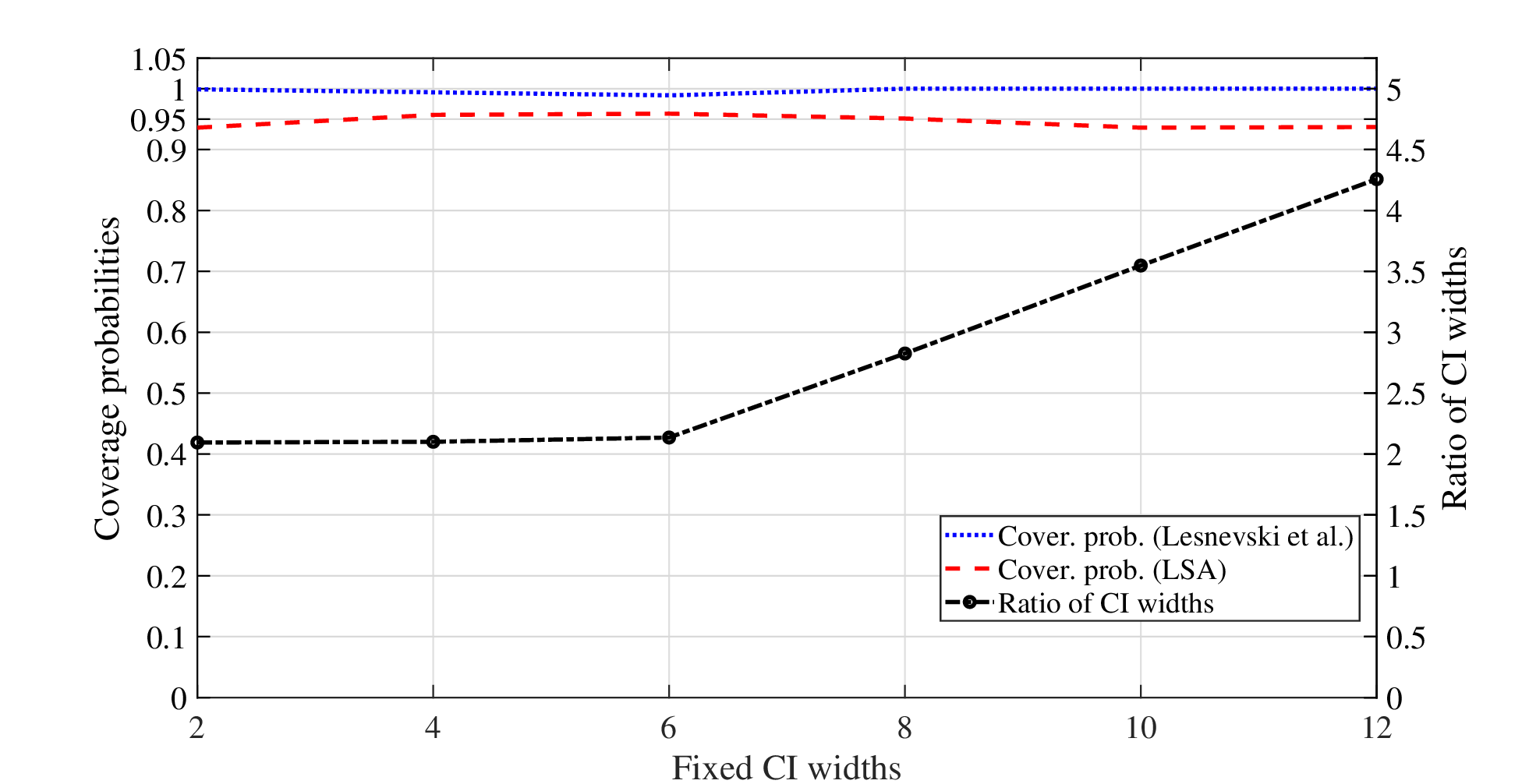}
		\caption{(Color online) Comparison of CIs}
		\label{fig:RiskMeasure}
	\end{center}
\end{figure}

In the comparison of CIs, we observed that the widths of the GA and the LSA CIs are quite close as they have similar variances. Hence, we only report the result of the LSA CIs when comparing their widths to the method of Lesnevski et al. (2007). The comparison results are summarized in Figure \ref{fig:RiskMeasure}, w.r.t. varying settings of the fixed width by Lesnevski et al. (2007). From the figure, it can be seen that the LSA CIs achieve the nominal coverage probability ($95\%$), while the CIs constructed by Lesnevski et al. (2007) tend to be conservative, yielding coverage probabilities close to $100\%$. The ratio between the width of the CI of Lesnevski et al. (2007) and that of the LSA CIs is also reported. It can be observed that the widths of the CIs of Lesnevski et al. (2007) are more than 2 times larger than LSA CIs, especially when the fixed width is set to be large. The ratio can be as large as 4 times when the fixed width is large, implying that the LSA method produces much narrower CIs and thus has better quality.

\subsection{Hypothesis Testing in Clinical Trials}\label{sub: ClinicalTrials}

In a clinical trial, suppose that one is interested in comparing $K-1$ treatments with unknown mean effects $\mu_1,\ldots, \mu_{K-1}$ to a control treatment with known mean effect $\mu_0$. To examine the performances of different hypothesis tests, we consider respectively the GA-based and LSA-based testing statistics in (\ref{eqn:TestGA}) and (\ref{eqn:TestLSA}), based on the proposed single hypothesis test in (\ref{eqn:Temp2Num}). We further consider a benchmarking test with multiple hypotheses in the form of (\ref{eqn:Temp1Num}). The benchmarking test uses a fixed randomization design that randomly assigns a patient to the treatments with equal probabilities, referred to as FR test.  

Our numerical setting is the same as Mozgunov and Jaki (2020), where the outcomes of treatments for patients follow Bernoulli distributions with success probabilities $(\mu_0,\mu_1,\ldots,\mu_{K-1})$, $K=4$, $n=423$, and overall type I error equal to $5\%$. In addition to determining whether some treatments outperform the control, another important goal in the clinical trial is to treat the patients as effectively as possible during the trial (Villar et al. 2015). To quantify to what extent this goal is achieved, two measures are commonly used, including the expected proportion of patients assigned to the best treatment (POB) and the expected number of patient success (ENS); see, e.g., Mozgunov and Jaki (2020).

To measure the performance of each test, we estimate its POB and ENS. Furthermore, for a case with true means $(\mu_0, \mu_1, \mu_2, \mu_3)=(0.31, 0.27, 0.28, 0.29)$ where $H_0$ is true, we estimate its overall type I error, while we estimate its power for another case with true means $(\mu_0, \mu_1, \mu_2, \mu_3)=(0.3, 0.3, 0.3, 0.5)$ where $H_0$ is false. Comparison results are summarized in Table \ref{tab:trials}, based on $10^4$ independent replications. 

It can be seen from Table \ref{tab:trials} that the GA-based and LSA-based tests have an advantage over the FR test in that they consistently assign more patients to the best treatment and have higher ENS, and are thus preferable when taking into account the well-being of the patients. From a statistical perspective, when $H_0$ is true (i.e., Case 1), the estimated type I error of the LSA-based test is close to the nominal one ($5\%$), while the rest tend to be conservative. Conservativeness of the GA-based test comes from its negative bias that is not negligible with a sample size of $423$, and that of the FR test is due to drawbacks of Bonferroni correction. In Case 2 when $H_0$ is false, both the GA-based and LSA-based tests have much higher estimated power compared to the FR test. In summary, the LSA-based test outperforms the rest of the tests significantly in all the considered dimensions including POB, ENS, type I error, and statistical power.  

\begin{table}[]%\footnotesize
\renewcommand\arraystretch{1.2}
\centering
\caption{Comparison of different hypothesis tests.}
\vspace{0.1in}
\label{tab:trials}
\begin{tabular}{cccccccc}
\hline
 %{\rule[-1mm]{0mm}{17mm}}
          & \multicolumn{3}{c}{Case 1: $(\mu_0, \mu_1, \mu_2, \mu_3)=(0.31, 0.27, 0.28, 0.29)$}&&\multicolumn{3}{c}{Case 2: $(\mu_0, \mu_1, \mu_2, \mu_3)=(0.3, 0.3, 0.3, 0.5)$}\\
             \cline{2-4}\cline{6-8}
Test                 &Est. Type I     &POB (s.e.)         &ENS (s.e.)       &&Est. power &POB (s.e.)    &ENS (s.e.)  \\
\hline
 FR               &0.014          &0.25(0.0001)       &121.6(0.03)   &&0.814    &0.25(0.0001)    &148.1(0.03) \\
 GA               &0.011          &0.34(0.0006)       &123.1(0.03)   &&0.966    &0.78(0.0006)    &192.7(0.02)    \\
 LSA              &0.044          &0.34(0.0006)       &123.1(0.03)   &&0.970    &0.78(0.0006)    &192.7(0.02)    \\
\hline
\end{tabular}
\end{table}

\subsection{Robust Analysis of Stochastic Simulation: A Call Center Case}\label{sec:RobustAnalysis}

An application of the proposed maximum mean estimator is on robust analysis of stochastic simulation. Robust analysis aims to quantify the impact of input uncertainty (i.e., misspecification of the input models) inherent in stochastic simulation. For this purpose, Ghosh and Lam (2019) and Blanchet and Murthy (2019) considered worst-case performance bounds, i.e., maximum (and/or minimum) of the expected simulation output taken over a set of plausible input models. While these authors consider a collection of plausible input models in the space of continuous distributions by imposing constraints on their moments, distances from a baseline model or some other summary statistics, we consider in this paper a discrete setting where a finite number of plausible input models are specified by fitting to input data.

Specifically, the expected simulation outputs corresponding to a set of $K$ plausible input models are denoted by $\mu_1, \ldots, \mu_K$, respectively. From each input model, a simulation output with noise can be generated. We are interested in estimating the following worst-case bounds,
\[\mu^{\dag} = \min_{k=1,\ldots, K}\mu_k,\quad \text{and} \quad \mu^*=\max_{k=1,\ldots,K}\mu_k,\]which help to quantify the input uncertainty of the plausible input models.

\begin{figure}[h!]
  \centering
  \includegraphics[width=1\textwidth]{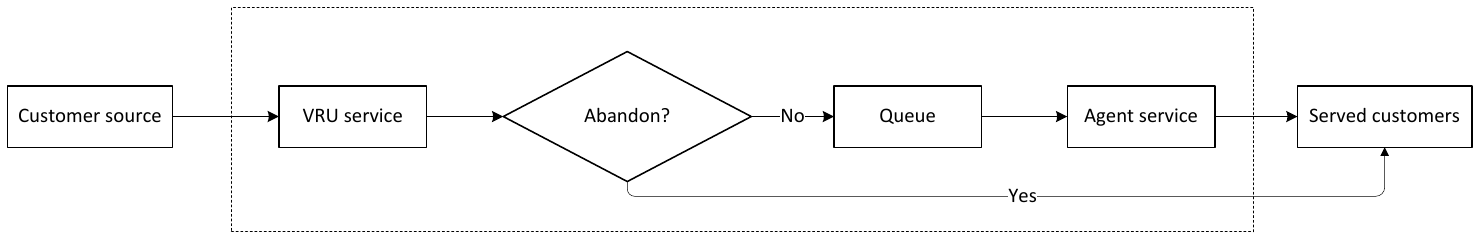}\\
  \caption{Logic flow of the call center of an anonymous bank}\label{fig: call center}
\end{figure}

The maximum mean estimators studied in this paper can be applied to estimate $\mu^*$, as well as $\mu^{\dag}$ by taking a negative sign on the simulation outputs. To illustrate how it works, we consider a case of a call center from an anonymous bank with real data available over a period of 12 months from January 1999 till December 1999 (Mandelbaum 2014). Workflows of the call center are described as follows. Customer phone calls arrive and ask for a banking service on their own. The call center consists of both automated service (computer-generated voice information, VRU) and agent service (AS) provided by eight agents. After the VRU service, customers either to a queue before AS, or directly to receive AS, or leave the system directly (see Figure \ref{fig: call center}). While privoviding VRU service seven days a week, 24 hours a day, the call center is staffed during 7:00 and 0:00 during weekdays. The information of each phone is stockpile, and each call has 17 fields in the data, e.g., Call ID (5 digits), VRU entry (6 digits), VRU exit (6 digits), Queue start (6 digits), Queue exit (6 digits), Outcome (AS or Hang), AS start (6 digits), AS exit (6 digits), etc.

We build up a simulation model with the logic flow described in Figure \ref{fig: call center}, to study customers' expected waiting time during the interval 10:00-12:00 on a weekday. Key inputs to the simulation include interarrival (IA) distribution, and service time distributions for both VRU and AS services. During input modeling, we consider a base dataset and an extended dataset, including input data during the specified interval on January 1, 1999 (with 172 phone calls), and the combined data on both January 1 and 8 (with 311 phone calls in total), respectively. Comparison of simulation performances based on these two datasets helps to understand the impact of input data size on the magnitude of input uncertainty. Arena's Input Analyzer tool is used to fit the input distributions, and Table \ref{table: fitted distribution} presents three candidate distributions for each of IA and AS times, all with $p$-value larger than 0.01, while empirical distribution of VRU service time is used due to unsatisfactory $p$-values for parametric distributions given by Inout Analyzer. The abandon rates for the two scenarios are 7.1\% and 5.5\%, respectively, according to the proportion of calls that hang out after VRU service.

\begin{table}[h!]
\caption{Fitted distributions with $p$-value larger than 0.01 for IA and AS times.}
\label{table: fitted distribution}
\centering
\vspace{10pt}
\begin{tabular}{lll}
\hline
Dataset   & IA time                                 & AS time                                            \\ \hline
Base dataset & -0.001 + ERLA(41.9, 1)                  & -0.001 + ERLA(225, 1)                                        \\
           & -0.001 + EXPO(41.9)                     & -0.001 + EXPO(225)                                            \\
           & -0.001 + WEIB(38.4, 0.82)               & -0.001 + WEIB(225, 0.992)                                    \\ \hline
Extended dataset  & -0.001 + ERLA(41.8, 1)                  & -0.001 + ERLA(222, 1)                                            \\
           & -0.001 + EXPO(41.8)                     & -0.001 + EXPO(222)                                           \\
           & -0.001 + WEIB(40.9, 0.951)              & -0.001 + WEIB(222, 0.996)                                    \\ \hline
\end{tabular}%
\end{table}

For each dataset, we specify the plausible input models as all possible combinations of the three candidate distributions for IA and AS times, while empirical distribution is employed for the VRU service time. In total there are $K=9$ ($=3\times 3$) plausible models for each dataset, and the proposed GA and LSA estimators are applied to estimate $\mu^{\dag}$ and $\mu^*$, the worst-case bounds for expected waiting time. Using real data, we also estimate the true value of the expected waiting time.

\begin{table}[H]%\footnotesize
\renewcommand\arraystretch{1.2}
\centering
\caption{Comparison of 95\% CIs of minimum mean $\mu^{\dag}$, maximum mean $\mu^{*}$, the plausible performance range $[\widehat{\mu}^\dag, \widehat{\mu}^*]$, and the gap between $\mu^{*}$ and $\mu^{\dag}$ for the base and extended datasets.}
\vspace{0.1in}
\label{table: CI results}
\begin{threeparttable}
\begin{tabular}{cccccccc}
%{p{2.5cm}<{\centering}p{0.8cm}<{\centering}p{0.8cm}<{\centering}p{0.8cm}<{\centering}p{0.8cm}<{\centering}p{0.001cm}p{0.7cm}<{\centering}p{0.7cm}<{\centering}p{0.7cm}<{\centering}p{0.8cm}<{\centering}}
\hline
          && \multicolumn{2}{c}{LSA}&&\multicolumn{2}{c}{GA}&\multirow{2}{*}{\centering True value\tnote{1}}\\
             \cline{3-4}\cline{6-7}
Dataset    & $n$                        &$10^4$        &$10^5$        &&$10^4$        &$10^5$  &\\
\hline
\multirow{4}{1.3cm}{\centering Base dataset}
           & CI of $\mu^\dag$           &[9.99, 10.87] &[10.04, 10.24]&&[10.40, 10.94]&[10.20, 10.37]&10.13\\
           & CI of $\mu^*$              &[13.70, 14.90]&[14.71, 14.97]&&[13.89, 14.60]&[14.64, 14.87]&14.81\\
&$[\widehat{\mu}^\dag, \widehat{\mu}^*]$&[10.43, 14.30]&[10.14, 14.84]&&[10.67, 14.25]&[10.29, 14.75]&\textbf{11.98}\\
& $\widehat{\mu}^* - \widehat{\mu}^\dag$& 3.87         & 4.70         && 3.58         & 4.46         &\\
\hline
\multirow{4}{1.3cm}{\centering Extended dataset}
           & CI of $\mu^\dag$           &[8.78, 9.13]  &[8.76, 8.91]  &&[9.29, 9.78]  &[8.76, 8.91]  &8.84\\
           & CI of $\mu^*$              &[11.69, 12.51]&[11.98, 12.21]&&[11.28, 11.86]&[11.86, 12.06]&12.10\\
&$[\widehat{\mu}^\dag, \widehat{\mu}^*]$&[9.05, 12.10] &[8.83, 12.10] &&[9.53, 11.57] &[8.84, 11.96] &\textbf{12.05}\\
& $\widehat{\mu}^* - \widehat{\mu}^\dag$&  3.05        &  3.27        && 2.04         & 3.12         &\\
\hline
\end{tabular}
\begin{tablenotes}
\footnotesize
\item[1] In the last column, the first two numbers for each dataset are the (approximately) true values of $\mu^\dag$ and $\mu^*$ estimated based a large ($10^5$) sample size for each plausible input model, and the bold numbers represent the average waiting times obtained by taking averages from the real data.
\end{tablenotes}
\end{threeparttable}
\end{table}

Table \ref{table: CI results} reports the constructed 90\% CIs for $\mu^{\dag}$ and $\mu^*$, the plausible performance range $[\widehat{\mu}^{\dag}, \widehat{\mu}^*]$, as well as the width $\widehat{\mu}^* - \widehat{\mu}^{\dag}$ obtained by the GA and LSA estimators. Based on the table, we highlight our observations as follows.
\begin{itemize}
\item Both the GA and the LSA estimators produce accurate estimates of $\mu^{\dag}$ and $\mu^*$ when sample size $n$ is sufficiently large, exemplified by CIs with narrow half-widths when $n=10^5$. However, the GA estimator may have larger bias when sample size is not sufficiently large. For instance, when $n=10^4$, the GA CIs of $\mu^{\dag}$ and $\mu^*$ lie above and below their true values (e.g. 10.13 and 14.81 for the base dataset) respectively, suggesting positive bias and negative bias for $\mu^{\dag}$ and $\mu^*$, respectively. By contrast, the LSA CIs have better performance as they cover the true values for both the base and extended datasets.

\item The estimated plausible ranges of $[\mu^{\dag}, \mu^*]$ resulting from the LSA estimator cover the true waiting times for both datasets, while those from the GA estimator cover true waiting times for $n=10^5$ but not for $n=10^4$ due to its larger bias. Compared to the GA estimator, the LSA estimator produces more accurate point estimates of both $\mu^{\dag}$ and $\mu^*$ (and thus the gaps of the worst-case bounds $\mu^*-\mu^{\dag}$) for all sample sizes and both datasets. Comparing the gaps of the worst-case bounds for these two datasets (with 172 and 311 pone calls respectively), it can be seen that increasing the size of the input dataset helps reducing significantly this gap and thus the input uncertainty, e.g., the gap is reduced from 4.70 to 3.27 using the LSA estimator with $n=10^5$.

\end{itemize}

\section{Conclusions}\label{sec:Conclusions}

We have studied two estimators, namely, the GA and LSA estimators, for estimating the maximum mean based on the UCB sampling policy. We have established strong consistency, asymptotic MSEs, and CLTs for both estimators, and constructed asymptotically valid CIs. Furthermore, a single hypothesis test has been proposed to solve a multiple comparison problem with application to clinical trials, and shown to significantly outperform an existing benchmarking test with multiple hypotheses, as the former exhibits lower type I error and higher statistical power. 

The LSA estimator is preferable over the GA estimator, mainly due to the fact that its bias decays at a rate much faster than that of the GA estimator. The superiority of the LSA estimator over the GA estimator and several competing estimators has been demonstrated through numerical experiments.

\section*{Acknowledgements}

This research is partially supported by the National Natural Science Foundation of China and  Research Grants Council of Hong Kong under Joint Research Scheme project N\_CityU 105/21 and GRF 11508217 and 11508620, and InnoHK initiative, The Government of the HKSAR, and Laboratory for AI-Powered Financial Technologies.

\newpage

\appendix
\renewcommand{\appendixname}{Appendix~\Alph{section}}
\setcounter{definition}{0}   %ÅœÃÃÃ£Â¿ÂªÃÅÂ±Ã ÂºÃ
\renewcommand{\thedefinition}{A\arabic{definition}}
\setcounter{lemma}{0}   %ÅœÃÃÃ£Â¿ÂªÃÅÂ±Ã ÂºÃ
\renewcommand{\thelemma}{A\arabic{lemma}}
%\setcounter{corollary}{0}   %ÅœÃÃÃ£Â¿ÂªÃÅÂ±Ã ÂºÃ
%\renewcommand{\thecorollary}{A\arabic{corollary}}

%\newpage

\section{Appendix}\label{sec:app}

\subsection{Lemmas}

\begin{lemma}\label{lem:Hoeffdingsub-Gaussian}
(Hoeffding bound, Proposition 2.5 in Wainwright, 2019) Suppose that random variables $\xi_i$, $i=1,...,n$, are independent, and $\xi_i$ has mean $\mu_i$ and sub-Gaussian parameter $\gamma_i$. Then for all $t\geq0$, we have
\begin{align*}
\P\left( \sum_{i=1}^n (\xi_i-\mu_i) \geq t \right) \leq \exp\left\{ -\frac{t^2}{2\sum_{i=1}^n \gamma_i^2}\right\}.
\end{align*}
\end{lemma}

\begin{lemma}\label{lem:SumBound}
Suppose that $X_i$'s are i.i.d. mean-zero random variables, and $N$ is a stopping time taking positive integer values with $N\le n$ w.p.1. Let $S_t = \sum_{i=1}^t X_i$. For any positive integer $p\ge 2$, if $\E|X_1|^p<\infty$, then
\[\E\left|S_N^p\right| = O(n^{p/2}).\] 
\end{lemma}

\begin{proof}[Proof of Lemma \ref{lem:SumBound}] First note that it suffices to show that $\left|\E(S_N^p) \right|= O(n^{p/2})$ for any $p\ge 2$, because these two are equivalent when $p$ is an even number, and if $p$ is an odd number, $\E\left|S_N^p\right| \le \left[\E\left( S_N^{p+1} \right)\right]^{p\over p+1} = O(n^{{p+1\over 2}\cdot {p\over p+1}}) = O(n^{p/2})$, following H\"older's inequality.   

For notational ease, denote $C_{q}^k = \binom{q}{k}$ and $\kappa_q = \E[X_1^q]$. Note that $\kappa_1 = 0$. Further let ${\cal F}_t$ denote the filtration up to $t$. 

We first note the following fact. With $c$ being a constant, $s$ a nonnegative integer, and $V_t$ a martingale, if a process $\{M_t, t\ge 1\}$ satisfies 
\[\E\left[ M_{t+1}|{\cal F}_t \right] = M_t + c t^s V_t,\]then it can be seen that 
\[W_t \triangleq M_t + \ba^\top \bt^{s+1} V_t\] 
is a martingale, by verifying the property that $\E[W_{t+1}|{\cal F}_t] = W_t$, where the shorthand notation $\bt^{s+1} \triangleq (t,\ldots,t^{s+1})^\top$, and the constant vector $\ba = (a_1,\ldots,a_{s+1})^\top$ satisfies 
\begin{equation}\label{eqn:TempSol}
a_{s+1} = -c/C_{s+1}^1,\quad \mbox{and}\quad a_{s-i} = -\sum_{k=2}^{i+2} a_{s-i-1+k}C_{s-i-1+k}^k/ C_{s-i}^1,\quad i = 0,\ldots, s -1.
\end{equation}

Define $V_{t,0}\equiv 1$ and $V_{t,1}=S_t$. For any nonnegative integer $q\ge 0$, if $\E\left[ S_{t+1}^q |{\cal F}_t \right] $ can be written in the form of
\begin{eqnarray}
\E\left[ S_{t+1}^q |{\cal F}_t \right] = S_t^q + \sum_{j\in {\cal S}_q} c(j,q) t^{s(j,q)}V_{t, r(j,q)},\label{eqn:Form}
\end{eqnarray}
where $c(j,q)$'s are constants,  $s(j,q)$ and $r(j,q)$ are nonnegative integers depending on $(j,q)$, and ${\cal S}_p$ is an index set with a finite number of elements, 
then, we define
\[V_{t,q} = S_t^q + \sum_{j\in {\cal S}_q} \ba(j,q)^\top \bt^{s(j,q)+1}V_{t, r(j,q)},\]where $\ba(j,q)$ is a contant vector obtained via (\ref{eqn:TempSol}) by setting $s = s(j,q)$ and $c=c(j,q)$. 

Then, using the aforementioned fact, the process $\{V_{t,q}, t\ge 1\}$ is a martingale, and therefore, 
\begin{equation}\label{eqn:LinearV}
S_t^q = V_{t,q} - \sum_{j\in {\cal S}_q} \ba(j,p)^\top \bt^{s(j,q)+1}V_{t, r(j,q)}
\end{equation}
is a linear combination of martingales. 

In the rest of this proof, we use an induction argument. Suppose that for any integer $q$ satisfying $2\le q\le p$, (a) and (b) hold, where 
\begin{enumerate}[label=\alph*)]
\item $\E\left[ S_{t+1}^q |{\cal F}_t \right] $ is of the form in (\ref{eqn:Form}), satisfying $V_{0,r(j,q)}=0$ when $r(j,q)\ge 1$, and $r(j,q)\le q-2$ for any $j\in {\cal S}_q$; and
\item $\E\left| N^{s(j,q)+1}V_{N, r(j,q)} \right| = O(n^{q/2})$ for any $s(j,q)$ and $r(j,q)$.,
\end{enumerate} 
we want to show that both (a) and (b) hold for $q=p+1$.

It should be pointed out that combining (\ref{eqn:LinearV}) and (b) yields $\left| \E(S_N^q) \right| = O(n^{q/2})$, because $\E[V_{N,q}] = V_{0,q} = 0$ due to the martingale stopping theorem, and
\begin{eqnarray*}
\E(S_N^q) &=&\E[V_{N,q}] - \sum_{j\in {\cal S}_q} \ba(j,p)^\top \E\left[\bN^{s(j,q)+1}V_{N, r(j,q)}\right]= O(n^{q/2}).
\end{eqnarray*}
Moreover, when $q=2$, it can be easily checked that $\E[S_{t+1}^2|{\cal F}_t] = S_t^2 + \kappa_2 V_{t,0}$ and $V_{t,2} = S_t^2 -\kappa_2 t V_{t,0}$, and thus both (a) and (b) are satisfied. 

Therefore, it suffices to complete the induction argument, i.e., (a) and (b) must hold for $q=p+1$ if they hold for any $2\le q\le p$. To see this, note that
\begin{eqnarray}\label{eqn:p_plus1}
\E\left[ S_{t+1}^{p+1} |{\cal F}_t\right] = \E\left[ \left( S_t + X_{t+1} \right)^{p+1} |{\cal F}_t\right]= S_t^{p+1} + \sum_{j=2}^{p+1} C_{p+1}^j \kappa_j S_t^{p+1-j}.
\end{eqnarray}

For any $j=2,\ldots,p-1$, combining (a) and (\ref{eqn:LinearV}) suggests that $S_t^{p+1-j}$ is a linear combination of $V_{t,p+1-j}$ and $V_{t, r(j,p+1-j)}$'s. Moreover, when $j=p$ and $p+1$, $S_t^{p+1-j}$ equals $V_{t,1}$ and $V_{t,0}$, respectively. Therefore, by (\ref{eqn:p_plus1}), assembling terms and rearranging the index set ${\cal S}_{p+1}$ lead to 
\begin{eqnarray*}
\E\left[ S_{t+1}^{p+1} |{\cal F}_t \right] = S_t^{p+1} + \sum_{j\in {\cal S}_{p+1}} c(j,p+1) t^{s(j,p+1)}V_{t, r(j,p+1)}
\end{eqnarray*}
for some $c(j,p+1)$, $s(j,p+1)$ and $r(j,p+1)$, where it can be verified that $V_{0, r(j,p+1)} = 0$ when $r(j,p+1)\ge 1$.

It can also be checked that
\begin{enumerate}[label=(\roman*)]
\item $r(j,p+1)=r(j_1,q)$ and $s(j,p+1) \le s(j_2,q)+1$ for some $q\le p-2$ and $j_1,j_2\in {\cal S}_q$; or,
\item $r(j,p+1)= q$ and $s(j,p+1)=0$ for some $q\le p-1$. 
\end{enumerate}

Because $r(j_1,q)\le q-2$, it is clear that $r(j,p+1)\le p-1$. Hence, (a) holds for $q=p+1$. 

Furthermore, define
\[V_{t,p+1} = S_t^{p+1} + \sum_{j\in {\cal S}_{p+1}} \ba(j,p+1)^\top \bt^{s(j,p+1)+1}V_{t, r(j,p+1)}.\]Then, $\{V_{t,p+1}, t\ge 1\}$ is a martingale. Using the martingale stopping theorem, we have $\E\left[ V_{N,p+1}\right]=V_{0,p+1}=0$, i.e.,
\begin{equation}\label{eqn:TempP}
\E\left[ S_N^{p+1} \right] =  -\sum_{j\in {\cal S}_{p+1}} \E\left[ \ba(j,p+1)^\top \bN^{s(j,p+1)+1}V_{N, r(j,p+1)}\right],
\end{equation}
where $\bN^{s(j,p+1)+1} = (N,\ldots,N^{s(j,p+1)+1})^\top$. 

When (i) is true,
\begin{eqnarray*}
\E\left| N^{s(j,p+1)+1}V_{N, r(j,p+1)} \right|\le \E\left| N^{s(j_2,q)+1}V_{N, r(j_1,q)} \right|\le n \E\left| N^{s(j_2,q)}V_{N, r(j_1,q)}\right| = O(n\cdot n^{q\over 2}) = O(n^{p\over 2}),
\end{eqnarray*}
where the last equality follows from the fact that $q\le p-2$.

When (ii) is true,  for some $q\le p-1$,
\begin{eqnarray*}
\E\left| N^{s(j,p+1)+1}V_{N, r(j,p+1)} \right| = \E\left| NV_{N, q} \right|.
\end{eqnarray*}
Note that $\left|V_{N,q}\right| \le |S_N^q| + \sum_{j\in {\cal S}_q} \left| \ba(j,q)^\top \bN^{s(j,q)+1}V_{N, r(j,q)}\right|$. If $q$ is an even number, then 
\[\E\left|V_{N,q}\right| \le \E(S_N^q) + \sum_{j\in {\cal S}_q} \E\left| \ba(j,q)^\top \bN^{s(j,q)+1}V_{N, r(j,q)}\right| = O(n^{q\over 2}) = O(n^{p-1\over 2}),\]where the first equality follows from (\ref{eqn:LinearV}) and (b). 
Otherwise, $q$ is an odd number, and by Cauchy-Schwartz inequality, 
\begin{eqnarray*}
\E\left|V_{N,q}\right| &\le& \E\left| S_N^q \right| + \sum_{j\in {\cal S}_q} \E\left| \ba(j,q)^\top \bN^{s(j,q)+1}V_{N, r(j,q)}\right|\\
&\le & \left[\E\left( S_N^{q+1} \right)\right]^{q\over q+1} + \sum_{j\in {\cal S}_q} \E\left| \ba(j,q)^\top \bN^{s(j,q)+1}V_{N, r(j,q)}\right|\\
&\le &O(n^{{q+1\over 2}\cdot{q\over q+1}}) + O(n^{q\over 2}) = O(n^{p-1\over 2}).
\end{eqnarray*}

Therefore no matter $q$ is even or odd, $\E\left| N^{s(j,p+1)+1}V_{N, r(j,p+1)} \right| = \E\left| NV_{N, q} \right|\le n\E\left|V_{N, q} \right|=O(n^{p+1\over 2})$. In other words, (b) holds for $q=p+1$, which completes the induction. 
\end{proof}

\subsection{Proof of Proposition \ref{prop:mom of Tk}}

To prove the proposition, we introduce a lemma.
\begin{lemma}\label{lem:IeqSum}
Given a real $\lambda\notin(0,1)$, for any positive integer $n$,
\begin{eqnarray*}
\sum_{j=1}^n j^\lambda  \begin{cases}
\displaystyle = n &if\ \lambda=0;\\
\displaystyle = \frac{n(n+1)}{2}  &if\ \lambda=1;\\
\displaystyle < \log\left(n+\frac{1}{2}\right)+\log2 &if\ \lambda=-1;\\
\displaystyle < \dfrac{1}{\lambda+1}\left[\left(n+\frac{1}{2}\right)^{\lambda+1}-\frac{1}{2^{\lambda+1}}\right] & if\ \lambda\neq-1.
\end{cases}
\end{eqnarray*}
\end{lemma}

\begin{proof}[Proof of Lemma \ref{lem:IeqSum}] The conclusion of the lemma for $\lambda=0$ or 1 is trivial, and we thus mainly examine the case when $\lambda\notin[0,1]$. Because $x^\lambda$ ($x>0$) is convex for $\lambda\notin[0,1]$, by Jensen's Inequality,
\begin{eqnarray*}
\int_{n-\frac{1}{2}}^{n+\frac{1}{2}}x^\lambda {\rm d}x>\left(\int_{n-\frac{1}{2}}^{n+\frac{1}{2}}x {\rm d}x\right)^\lambda=n^\lambda,\ {\rm if}\ \lambda\notin[0,1].
\end{eqnarray*}
It then follows that
\begin{eqnarray*}
\left(n+\frac{1}{2}\right)^{\lambda+1}-\left(n-\frac{1}{2}\right)^{\lambda+1}=(\lambda+1)\int_{n-\frac{1}{2}}^{n+\frac{1}{2}}x^\lambda {\rm d}x
\begin{cases}
 \displaystyle  > (\lambda+1)n^\lambda &{\rm if}\ -1<\lambda<0,\ \  {\rm or}  \ \ \lambda>1,\\
 \displaystyle  < (\lambda+1)n^\lambda,  &{\rm if}\ \lambda<-1,
\end{cases}
\end{eqnarray*}
and
\begin{eqnarray*}
\log\left(n+\frac{1}{2}\right)-\log\left(n-\frac{1}{2}\right)=\int_{n-\frac{1}{2}}^{n+\frac{1}{2}}x^{-1} {\rm d}x > n^{-1},\ {\rm if}\  \lambda=-1.
\end{eqnarray*}

Therefore, for $\lambda\notin[0,1]$ and $\lambda\neq-1$,
\begin{eqnarray*}
\sum_{j=1}^n j^\lambda <\dfrac{1}{\lambda+1}\sum_{j=1}^n \left[\left(j+\frac{1}{2}\right)^{\lambda+1}-\left(j-\frac{1}{2}\right)^{\lambda+1}\right]=\dfrac{1}{\lambda+1}\left[\left(n+\frac{1}{2}\right)^{\lambda+1}-\left(\frac{1}{2}\right)^{\lambda+1}\right],
\end{eqnarray*}
and for $\lambda=-1$,
\begin{eqnarray*}
\sum_{j=1}^n j^\lambda< \sum_{j=1}^n \left[\log\left(j+\frac{1}{2}\right)-\log\left(j-\frac{1}{2}\right)\right]=\log\left(n+\frac{1}{2}\right)-\log\left(\frac{1}{2}\right).
\end{eqnarray*}
\end{proof}

\noindent {\bf Proof of Proposition \ref{prop:mom of Tk}}.
We first the result with $p=1$ and then $p\geq2$.

(i) When $p=1$, we want to prove that for $k\neq k^*$
\begin{eqnarray}\label{eqn:Tk1}
\E T_k(n) \leq \frac{8\nu_n}{\Delta_k^2} + 5
\end{eqnarray}

Let
\begin{eqnarray*}
j^* = \max\left\{ j: T_k(j-1)\leq\frac{8\nu_n}{\Delta_k^2} \right\}.
\end{eqnarray*}
If $j^*=n$, then (\ref{eqn:Tk1}) is true. Therefore, we just have to consider $j^*<n$. Note that for $j>j^*$, we have
\begin{eqnarray}\label{eqn:Tk10}
T_k(j-1) > \frac{8\nu_n}{\Delta_k^2}.
\end{eqnarray}

For $j>j^*$ and $I_j=k$, we claim that at least one of the two following inequalities holds:
\begin{eqnarray}
\bar X_{k^*}[T_{k^*}(j-1)] + c_{j,T_{k^*}(j-1)} &\leq& \mu_{k^*},\label{eqn:Tk11}\\
\bar X_k[T_k(j-1)] - c_{j,T_k(j-1)} &>& \mu_k,\label{eqn:Tk12}
\end{eqnarray}
where $c_{j,s} = \sqrt{2\nu_j/s}$. In fact, by contradiction, if both the inequalities (\ref{eqn:Tk11}) and (\ref{eqn:Tk12}) are false when $j>j^*$ and $I_j=k$, we have
\begin{align*}
\bar X_{k^*}[T_{k^*}(j-1)] + c_{j,T_{k^*}(j-1)} >& \mu_{k^*} = \mu_k + \Delta_k > \mu_k + 2\sqrt{\frac{2\nu_n}{T_k(j-1)}}\\ \geq& \mu_k + 2\sqrt{\frac{2\nu_j}{T_k(j-1)}} \geq \bar X_k[T_k(j-1)] + c_{j,T_k(j-1)},
\end{align*}
where the second inequality is due to (\ref{eqn:Tk10}), which contradicts with $I_j=k$.

In the following, we prove (\ref{eqn:Tk1}). Note that
\begin{align}\label{eqn:Tk13}
&T_k(n) = T_k(j^*) + \sum_{j=j^*+1}^n \1\{I_j = k\} = T_k(j^*) + \sum_{j=j^*+1}^n \1\{I_j = k, \mbox{(\ref{eqn:Tk11}) or (\ref{eqn:Tk12}) true}\} \nonumber\\
&\leq T_k(j^*-1) + 1 + \sum_{j=j^*+1}^n \1\{\mbox{(\ref{eqn:Tk11}) or (\ref{eqn:Tk12}) true}\}
\leq \frac{8\nu_n}{\Delta_k^2} + 1 + \sum_{j=j^*+1}^n \1\{\mbox{(\ref{eqn:Tk11}) or (\ref{eqn:Tk12}) true}\}.
\end{align}
Taking expectation, it follows that
\begin{align*}
\E T_k(n) \leq& \frac{8\nu_n}{\Delta_k^2} + 1 + \sum_{j=j^*+1}^n \P \left(\mbox{(\ref{eqn:Tk11}) or (\ref{eqn:Tk12}) true} \right)\\
\leq& \frac{8\nu_n}{\Delta_k^2} + 1 + \sum_{j=j^*+1}^n \left[ \P \left(\mbox{(\ref{eqn:Tk11}) true} \right) + \P \left(\mbox{(\ref{eqn:Tk12}) true} \right) \right].
\end{align*}
Furthermore, by Lemma \ref{lem:Hoeffdingsub-Gaussian},
\begin{align}\label{eqn:Tk131}
\P \left(\mbox{(\ref{eqn:Tk11}) true} \right) =& \P \left(\bar X_{k^*}[T_{k^*}(j-1)] + c_{j,T_{k^*}(j-1)} \leq \mu_{k^*} \right)
\leq \P \left(\min_{s=1,...,j-1}\bar X_{k^*}[s] + c_{j,s} \leq \mu_{k^*} \right) \nonumber\\
=& \P \left(\cup_{s=1}^{j-1}\bar X_{k^*}[s] + c_{j,s} \leq \mu_{k^*} \right)
\leq \sum_{s=1}^{j-1}\P \left(\bar X_{k^*}[s] + c_{j,s} \leq \mu_{k^*} \right) \nonumber\\
\leq& \sum_{s=1}^{j-1} \exp\left\{-\frac{(s c_{j,s})^2}{2s \gamma_{k^*}^2} \right\}
= (j-1) \exp\left\{-\frac{\nu_j}{\gamma_{k^*}^2} \right\} < j^{-3},
\end{align}
where the last inequality is due to the condition $\nu_n\geq4\bar\gamma^2 \log n$. Similarly, $\P \left(\mbox{(\ref{eqn:Tk12}) true} \right)$ has the same bound. Therefore,
\begin{align*}
\E T_k(n) \leq \frac{8\nu_n}{\Delta_k^2} + 1 + 2 \sum_{j=j^*+1}^n j^{-3}
\leq \frac{8\nu_n}{\Delta_k^2} + 1 + 4 - \left(n+\frac{1}{2}\right)^{-2}
\leq \frac{8\nu_n}{\Delta_k^2} + 5,
\end{align*}
where the second inequality is by Lamma \ref{lem:IeqSum}.

(ii) Next, we prove the result with $p\geq2$. Denote
\begin{eqnarray}\label{def Rj}
R_j=\1\left\{I_j=k\right\}, j=j^*+1,...,n.
\end{eqnarray}
Then, by (\ref{eqn:Tk13}), we have
\begin{eqnarray}\label{Ieq:Tkn}
T_k(n)\leq \frac{8\nu_n}{\Delta_k^2} + 1 +\sum_{j=j^*+1}^n R_j.
\end{eqnarray}

For any integer $p\geq2$, by Minkowski's inequality, we have
\begin{eqnarray}\label{Ieq:Minkv}
\E T_k(n)^p&\leq& \E\left(\frac{8\nu_n}{\Delta_k^2} + 1 + \sum_{j=j^*+1}^n R_j\right)^p \leq \left\{\frac{8\nu_n}{\Delta_k^2} + 1 + \left[\E\left(\sum_{j=j^*+1}^n R_j\right)^p\right]^\frac{1}{p}\right\}^p.
\end{eqnarray}

In what follows, we analyze $\E\left(\sum_{j=j^*+1}^n R_j\right)^p$. Note that
\begin{eqnarray}\label{sum Rj p}
\lefteqn{\left(\sum_{j=j^*+1}^n R_j\right)^p = \sum_{j_1,...,j_p\in\{j^*+1,...,n\}} R_{j_1}\cdots R_{j_p}}\nonumber\\
&=& p!\sum_{j_1>\cdots>j_p} R_{j_1} \cdots R_{j_p} + \sum_{m=1}^{p-1}   \sum_{\begin{subarray}{l} i_1+\cdots+i_m=p\\ i_1,...,i_m\geq1\end{subarray}}\binom{p}{i_1,...,i_m} \sum_{j_1>\cdots>j_m}R_{j_1}^{i_1} R_{j_2}^{i_2}\cdots R_{j_m}^{i_m}\nonumber\\
&=& p!\sum_{j_1>\cdots>j_p} R_{j_1} \cdots R_{j_p} + \sum_{m=1}^{p-1}   \sum_{\begin{subarray}{l} i_1+\cdots+i_m=p\\ i_1,...,i_m\geq1\end{subarray}}\binom{p}{i_1,...,i_m} \sum_{j_1>\cdots>j_m}R_{j_1} R_{j_2}\cdots R_{j_m}
\end{eqnarray}
where $\binom{p}{i_1,...,i_m}=p!/(i_1!\cdots i_m!)$, and the last equality follows from the definition of $R_j$ in (\ref{def Rj}).

Consider $\sum_{j_1>\cdots>j_m} R_{j_1} \cdots R_{j_m}$ when $2\leq m\leq p$. By Lemma \ref{lem:IeqSum}, for a positive integer $\lambda$,
\begin{eqnarray}\label{IeqSum1}
\sum_{j=1}^n j^\lambda< \frac{1}{\lambda+1}\left(n+1\right)^{\lambda+1}.
\end{eqnarray}
Therefore, for $2\leq m\leq p$, we have
\begin{eqnarray*}
\lefteqn{\sum_{j_1>\cdots>j_m} R_{j_1} \cdots R_{j_m}=\sum_{j_1>\cdots>j_{m-1}} R_{j_1} \cdots R_{j_{m-1}}\sum_{j_m=j^*+1}^{j_{m-1}-1} R_{j_m}}\\
&\leq& \sum_{j_1>\cdots>j_{m-1}} R_{j_1} \cdots R_{j_{m-1}}\cdot j_{m-1}=\sum_{j_1>\cdots>j_{m-2}} R_{j_1} \cdots R_{j_{m-2}}\sum_{j_{m-1}=j^*+1}^{j_{m-2}-1} R_{j_{m-1}}\cdot j_{m-1}\\
&\leq&\sum_{j_1>\cdots>j_{m-2}} R_{j_1} \cdots R_{j_{m-2}}\sum_{j_{m-1}=j^*+1}^{j_{m-2}-1} j_{m-1}\leq\frac{1}{2}\sum_{j_1>\cdots>j_{m-2}} R_{j_1} \cdots R_{j_{m-2}} \cdot j_{m-2}^2,
\end{eqnarray*}
where the first inequality follows from the fact that $|R_j|\le 1$, and the last inequality follows from (\ref{IeqSum1}) by letting $\lambda=1$ and $n=j_{m-2}-1$.

Iteratively applying the above argument leads to
\begin{eqnarray}\label{ineq Rj}
\sum_{j_1>\cdots>j_m} R_{j_1} \cdots R_{j_m}\leq\frac{1}{(m-1)!}\sum_{j_1=j^*+1}^n R_{j_1}\cdot j_1^{m-1}.
\end{eqnarray}
By (\ref{eqn:Tk13}) and (\ref{eqn:Tk131}), we know that
\begin{align*}
\E R_{j_1} \leq \P \left(\mbox{(\ref{eqn:Tk11}) true} \right) + \P \left(\mbox{(\ref{eqn:Tk12}) true} \right) \leq 2j_1^{-3}.
\end{align*}

Hence, combining with (\ref{ineq Rj}), we have, for $2\leq m\leq p$,
\begin{eqnarray*}
\E\sum_{j_1>\cdots>j_m} R_{j_1} \cdots R_{j_m}\leq\frac{1}{(m-1)!}\sum_{j_1=j^*+1}^n j_1^{m-1}\E R_{j_1}\leq\frac{2}{(m-1)!}\sum_{j_1=1}^n j_1^{m-4}.
\end{eqnarray*}
By Lemma \ref{lem:IeqSum}, we have
\begin{eqnarray}\label{Ieq:Exp Rjm}
\E\sum_{j_1>\cdots>j_m} R_{j_1} \cdots R_{j_m} \leq \begin{cases}
\displaystyle 2\left(2-\frac{1}{n+\frac{1}{2}}\right) &if\ m=2;\\
\displaystyle \log\left(n+\frac{1}{2}\right)+\log2 &if\ m=3;\\
\displaystyle \frac{2}{(m-1)!}\frac{1}{m-3} \left(n+\frac{1}{2}\right)^{m-3} & if\ m\geq4.\\
\end{cases}
\end{eqnarray}

By (\ref{sum Rj p}) and (\ref{Ieq:Exp Rjm}), it follows that for $p\geq4$,
\begin{align*}
\E\left(\sum_{j=j^*+1}^n R_j\right)^p&= p!\cdot\E\sum_{j_1>\cdots>j_p} R_{j_1} \cdots R_{j_p} + \sum_{m=1}^{p-1}   \sum_{\begin{subarray}{l} i_1+\cdots+i_m=p\\ i_1,...,i_m\geq1\end{subarray}}\binom{p}{i_1,...,i_m} \cdot\E\sum_{j_1>\cdots>j_m} R_{j_1} R_{j_2}\cdots R_{j_m}\\
&\leq \frac{2p}{p-3} (n+1)^{p-3}+O\left((n+1)^{p-4}\right).
\end{align*}
Plugging the above into (\ref{Ieq:Minkv}), we have, for $p\geq4$,
\begin{eqnarray*}
\E T_k(n)^p\leq \left\{\frac{8\nu_n}{\Delta_k^2} + 1 + \left[\frac{2p}{p-3} (n+1)^{p-3}+O\left(n^{p-4}\right)\right]^\frac{1}{p} \right\}^p.
\end{eqnarray*}

Similarly, for $p=3$, we have
\begin{align*}
\E\left(\sum_{j=j^*+1}^n R_j\right)^3&= 6\left[ \log\left(n+\frac{1}{2}\right)+\log2 \right] + O(1)
\leq 6\log\left(n+1\right) + O(1).
\end{align*}
Hence,
\begin{eqnarray*}
\E T_k(n)^3 \leq \left\{\frac{8\nu_n}{\Delta_k^2} + 1 + \left[6\log\left(n+1\right) + O(1)\right]^\frac{1}{3} \right\}^3.
\end{eqnarray*}

For $p=2$, it follows from (\ref{IeqSum1}) and (\ref{Ieq:Exp Rjm}) that
\begin{align*}
&\E\left(\sum_{j=j^*+1}^n R_j\right)^2= \sum_{j=j^*+1}^n\E R_j^2 + 2\E\sum_{j_1>j_2} R_{j_1}R_{j_2}= \sum_{j=j^*+1}^n\E R_j + 2\E\sum_{j_1>j_2} R_{j_1}R_{j_2}\\
&\leq \sum_{j_1=j^*+1}^n 2 j_1^{-3} + 4\left[2-\frac{1}{n+\frac{1}{2}}\right]
< 4 - \left(n+\frac{1}{2}\right)^{-2} + 4\left[2-\frac{1}{n+\frac{1}{2}}\right]
< 12,
\end{align*}
and thus
\begin{eqnarray*}
\E T_k(n)^2\leq \left\{ \frac{8\nu_n}{\Delta_k^2} + 1 + \sqrt{12} \right\}^2 < \left\{ \frac{8\nu_n}{\Delta_k^2} + 5 \right\}^2.
\end{eqnarray*}

The proof is complete.

\subsection{Proof of Theorem \ref{thm:Tk a.s.}}\label{sec:ProofSamplingRatio}

When $k\neq k^*$, by Markov's inequality and Proposition \ref{prop:mom of Tk}, for any $\varepsilon>0$ and $p\geq3$,
\begin{eqnarray*}
\P\left( \frac{T_k(n)}{n}>\varepsilon\right)\leq \frac{1}{\varepsilon^p}\E\left[\frac{T_k(n)}{n}\right]^p
&\leq& \frac{1}{\varepsilon^p}\left\{ \frac{8}{\Delta_k^2} \frac{\nu_n}{n} + \frac{1}{n} + \left[ \frac{2p}{p-3}\frac{1}{n^3}+O\left(n^{-4}\right)\right]^\frac{1}{p} \right\}^p\\
&\leq&\frac{2^{p-1}}{\varepsilon^p}\left[ \left( \frac{8}{\Delta_k^2} \frac{\nu_n}{n} + \frac{1}{n}\right)^p + \frac{2p}{p-3}\frac{1}{n^3}+O\left(n^{-4}\right)\right].
\end{eqnarray*}

Because $\nu_n\leq n^{1-\delta}$, it follows that for $p\geq4$,
\begin{eqnarray}
\P\left( \frac{T_k(n)}{n}>\varepsilon\right)&\leq& \frac{2^{p-1}}{\varepsilon^p}\left[\left(\frac{8}{\Delta_k^2 n^\delta} + \frac{1}{n}\right)^p + \frac{2p}{p-3}\frac{1}{n^3}+O\left(n^{-4}\right)\right]\nonumber\\
&\leq& \frac{2^{p-1}}{\varepsilon^p}\left[2^{p-1}\left(\frac{8}{\Delta_k^2 n^{p\delta}} + \frac{1}{n^p}\right) + \frac{2p}{p-3}\frac{1}{n^3}+O\left(n^{-4}\right)\right].\label{eqn:Temp2}
\end{eqnarray}
Let $p>\max\{1/\delta,4\}$. Then, for any $\epsilon>0$,
\begin{eqnarray*}
\sum_{n=1}^\infty\P\left( \frac{T_k(n)}{n}>\varepsilon\right)< \infty.
\end{eqnarray*}
By Borel-Cantelli lemma, $T_k(n)/n\overset{a.s.}{\longrightarrow}0$ as $n\to\infty$, for $k\neq k^*$, and thus
\begin{eqnarray*}
\frac{T_{k^*}(n)}{n}=1-\sum_{k\neq k^*}\frac{T_k(n)}{n}\overset{a.s.}{\longrightarrow}1.
\end{eqnarray*}

\subsection{Proof of Theorem \ref{thm: Istar2kstar}}\label{sec:ProofInStar}

It follows from Theorem \ref{thm:Tk a.s.} that
\begin{eqnarray*}
\P\left(\omega\in\Omega: \lim_{n\to\infty}\frac{T_{k}(n)(\omega)}{n}=0,\ \forall k\neq k^*\right)=1.
\end{eqnarray*}
Because $\sum_{k=1}^K T_{k}(n)=n$, the event $\lim_{n\to\infty}\frac{T_{k}(n)(\omega)}{n}=0,\ \forall k\neq k^*$ implies that $\lim_{n\to\infty}\frac{T_{k^*}(n)}{n}=1$. Therefore,
\begin{eqnarray*}
\P\left(\omega\in\Omega: \lim_{n\to\infty}\frac{T_{k}(n)(\omega)}{n}=0,\ \forall k\neq k^*,\ {\rm and}\ \lim_{n\to\infty}\frac{T_{k^*}(n)(\omega)}{n}=1\right)=1.
\end{eqnarray*}
Then, $\P(A)=1$, where
\begin{eqnarray*}
A=\left\{\omega\in\Omega: \lim_{n\to\infty}\left[\frac{T_{k^*}(n)(\omega)}{n}-\frac{T_k(n)(\omega)}{n}\right]=1,\ \forall k\neq k^*\right\}.
\end{eqnarray*}

In other words, for any $\omega\in A$ and $0<\varepsilon<1$, there exists $N=N(\omega,\varepsilon)\in \N$, such that for any $n>N$,
\begin{eqnarray*}
0<1-\varepsilon<\frac{T_{k^*}(n)(\omega)}{n}-\frac{T_k(n)(\omega)}{n}<1+\varepsilon,
\end{eqnarray*}
which implies
\begin{eqnarray*}
\max\left\{{T_1(n)(\omega)},...,{T_K(n)(\omega)}\right\}={T_{k^*}(n)(\omega)}.
\end{eqnarray*}
By the uniqueness assumption of $k^*$, we have
\begin{eqnarray}\label{a.s. In0}
\P\left(\omega\in\Omega: \forall 0<\varepsilon<1, \exists N=N(\omega,\varepsilon)\in \N, s.t.\ \forall n>N, I^*_n(\omega)=k^*\right)=1,
\end{eqnarray}
i.e., $I^*_n\overset{a.s.}{\longrightarrow}k^*$ as $n\rightarrow\infty$.

\subsection{Proof of Theorem \ref{thm:MSE tilde M}}\label{sec:ProofMSE_GA}

We analyze the bias and variance of $\widetilde M_n$ separately.

\textbf{The Bias}. We assert that $T_k(n)$ is a stopping time. In fact, given the filtration
generated by $\{X_{ij},i\neq k, j = 1,2,\dots\}$, denoted by
$\sigma\left(X_{ij},i\neq k, j = 1,2,\dots\right)$, one has
\[\{T_k(n)\ge m\}=\{T_k(n)\le m-1\}^c\in
\sigma\left(X_{k1},\dots,X_{k,m-1}\right).\]Hence $T_k(n)$ is a
stopping time with respect to $\{X_{kj},j\ge 1\}$. Then by Wald's equation,
\begin{eqnarray*}
\E\left[\sum_{j=1}^{T_k(n)}X_{kj}\right]= \mu_k \E\left[T_k(n)\right],
\end{eqnarray*}
leading to
\[\E \widetilde M_n={1\over n}\sum_{k=1}^{K}\mu_k\E\left[T_k(n)\right].\]

Thus,
\begin{eqnarray*}
{\rm Bias}\left(\widetilde M_{n}\right)&=&\E \widetilde M_n - \mu^*=\left(-1+\frac{\E T_{k^*}(n)}{n}\right)\mu_{k^*}+ \sum_{k\neq k^*}\frac{\E T_k(n)}{n} \mu_k\\
&=&-\sum_{k\neq k^*}\frac{\E T_k(n)}{n}\mu_{k^*}+\sum_{k\neq k^*}\frac{\E T_k(n)}{n} \mu_k=\sum_{k\neq k^*}\frac{\E T_k(n)}{n} \left(\mu_k-\mu_{k^*}\right).
\end{eqnarray*}
Because $\mu_k<\mu_{k^*}$ for $k\neq k^*$, we have ${\rm Bias}\left(\widetilde M_{n}\right)<0$.

By Proposition \ref{prop:mom of Tk},
\begin{eqnarray}\label{bias tilde M}
\left|{\rm Bias}\left(\widetilde M_{I_n^*}\right)\right| \leq\sum_{k\neq k^*}\frac{1}{n}\left( \frac{8\nu_n}{\Delta_k^2} + 5\right) \left(\mu_{k^*} - \mu_k\right)\leq\sum_{k\neq k^*}\frac{1}{n}\left( \frac{8\nu_n}{\Delta_k^2}+5\Delta_k\right)=O\left(\frac{\nu_n}{n}\right).
\end{eqnarray}

\textbf{The Variance}. Note that
\begin{eqnarray}\label{eqn:Temp1}
&&\Var\left[\widetilde M_n\right]=\E\left(\widetilde M_n - \E\left[\widetilde M_n\right]\right)^2\nonumber\\
&=& \E\left\{\left(\widetilde M_n-{1\over n}\sum_{k=1}^{K}\mu_k T_k(n)\right)+\left({1\over n}\sum_{k=1}^{K}\mu_k T_k(n)-{1\over n}\sum_{k=1}^{K}\mu_k\E T_k(n)\right)\right\}^2\nonumber\\
&\leq & 2 \E\left(\widetilde M_n-{1\over n}\sum_{k=1}^{K}\mu_k T_k(n)\right)^2+2\E\left({1\over n}\sum_{k=1}^{K}\mu_k T_k(n)-{1\over n}\sum_{k=1}^{K}\mu_k\E T_k(n)\right)^2.
\end{eqnarray}

We analyze the two terms on the RHS of (\ref{eqn:Temp1}) separately. By the Wald's Lemma (see Theorem
13.2.14 in Athreya and Lahiri 2006),
\begin{eqnarray}\label{eqn:Wald lemma0}
\E\left(\sum_{j=1}^{T_k(n)}X_{kj}-\mu_kT_k(n)\right)^2=\sigma_k^2\E T_k(n).
\end{eqnarray}

Then, by (\ref{eqn:tilde M decomp}) and Cauchy-Schwarz inequality, it can be seen that
\begin{eqnarray}\label{eqn:Temp11}
\E\left(\widetilde M_n-{1\over n}\sum_{k=1}^{K}\mu_k T_k(n)\right)^2\leq \frac{K}{n^2} \sum_{k=1}^{K}\E\left(\sum_{j=1}^{T_k(n)}X_{kj}-\mu_kT_k(n)\right)^2=\frac{K}{n^2} \sum_{k=1}^{K}\sigma_k^2 \E T_k(n)\nonumber\\
=\frac{K}{n^2} \sigma_{k^*}^2 \E T_{k^*}(n) + \frac{K}{n^2}\sum_{k\neq k^*}\E T_k(n)\leq \frac{K}{n} \sigma_{k^*}^2 + \frac{K}{n^2}\sum_{k\neq k^*}\E T_k(n),
\end{eqnarray}
where the last inequality follows from the fact that $T_{k^*}(n)\leq n$.

We then analyze the second term on the RHS of (\ref{eqn:Temp1}).
Note that, by Cauchy-Schwarz inequality,
\begin{eqnarray}\label{eqn:Temp12}
&&\E\left({1\over n}\sum_{k=1}^{K}\mu_k T_k(n)-{1\over n}\sum_{k=1}^{K}\mu_k\E T_k(n)\right)^2= {1\over n^2}\E\left(\sum_{k=1}^{K}\mu_k T_k(n)-\sum_{k=1}^{K}\mu_k\E T_k(n)\right)^2\nonumber\\
&\leq & {K\over n^2}\sum_{k=1}^{K}\mu_k^2\E\left[T_k(n)-\E T_k(n)\right]^2=  {K\over n^2}\sum_{k=1}^{K}\mu_k^2 \Var\left[T_k(n) \right]\nonumber\\
&\leq & {K\over n^2}\sum_{k\neq k^*}\mu_k^2\Var\left[T_k(n)\right]+{K\over n^2}\mu_{k^*}^2 \Var\left[T_{k^*}(n)\right].
\end{eqnarray}
Furthermore,
\begin{eqnarray*}
\Var\left[T_{k^*}(n) \right] = \Var\left[n-\sum_{k\neq k^*}T_k(n)\right] = \Var\left[\sum_{k\neq k^*}T_k(n)\right] \leq \E\left[\sum_{k\neq k^*}T_k(n)\right]^2 \leq (K-1)\sum_{k\neq k^*}\E T_k^2(n).
\end{eqnarray*}
Therefore,
\begin{eqnarray}\label{eqn:Temp121}
{\rm RHS\ of\ }(\ref{eqn:Temp12})\leq{K\over n^2}\sum_{k\neq k^*}\mu_k^2\E T_k^2(n)+{K(K-1)\over n^2}\mu_{k^*}^2 \sum_{k\neq k^*}\E T_k^2(n)\leq \frac{K^2}{n^2}\bar\mu^2\sum_{k\neq k^*} \E T_k^2(n),
\end{eqnarray}
where $\bar\mu^2=\max\{\mu_1^2,...,\mu_K^2\}$.

Combining (\ref{eqn:Temp1}), (\ref{eqn:Temp11})-(\ref{eqn:Temp121}) and Proposition \ref{prop:mom of Tk} yields
\begin{align}\label{var tilde M}
\Var\left[\widetilde M_n\right]\leq& \frac{2K}{n} \sigma_{k^*}^2 + \frac{2K}{n^2}\sum_{k\neq k^*}\E T_k(n) + \frac{2K^2\bar\mu^2}{n^2}\sum_{k\neq k^*} \E T_k^2(n)\nonumber\\
\leq& \frac{2K}{n} \sigma_{k^*}^2 + \frac{2K}{n^2}\sum_{k\neq k^*}\left( \frac{8\nu_n}{\Delta_k^2} + 5\right) + \frac{2K^2\bar\mu^2}{n^2} \sum_{k\neq k^*} \left\{ \frac{8\nu_n}{\Delta_k^2} + 5 \right\}^2\nonumber\\
\leq& \frac{2K}{n} \sigma_{k^*}^2 + \frac{2K^2 c\nu_n}{n^2} + \frac{2K^3\bar\mu^2}{n^2} c^2 \nu_n^2
\leq \frac{2K}{n} \sigma_{k^*}^2 + \frac{2K^3(\bar\mu^2+1)}{n^2} c^2 \nu_n^2,
\end{align}
where the third inequality is because of the condition $\frac{8\nu_n}{\Delta_k^2} + 5\leq c\nu_n$ for $k\neq k^*$.

Incorporating the bias in (\ref{bias tilde M}), the variance in (\ref{var tilde M}), and the condition $\nu_n\leq n^{1/2-\delta}$ with $0<\delta<1/2$, we have
\begin{eqnarray*}
{\rm MSE}\left(\widetilde M_n\right)=\Var\left[\widetilde M_n\right]+{\rm Bias}\left(\widetilde M_{I_n^*}\right)^2\le  \frac{2K}{n} \sigma_{k^*}^2 + o\left(\frac{1}{n}\right).
\end{eqnarray*}

\subsection{Proof of Theorem \ref{thm:CLT tilde M}}\label{sec:ProofCLT_GA}

To prove the result, we first define a martingale difference array and introduce a lemma on asymptotic normality for martingale different arrays.

\begin{definition}\label{def:md}
{\rm (Definition 16.1.1 in Athreya and Lahiri (2006))}
Let $\left\{Y_n\right\}_{n\geq1}$ be a collection of random variables on a probability space $(\Omega,\cF, \P)$ and let $\left\{\cF_n\right\}_{n\geq1}$ be a filtration. Then,
$\left\{Y_n, \cF_n\right\}_{n\geq1}$ is called a martingale difference array if $Y_n$ is $\cF_n$-measurable and $\E[Y_n|\cF_{n-1}]=0$ for each $n\geq1$.
\end{definition}

\begin{lemma}[Theorem 16.1.1 of Athreya and Lahiri (2006)]\label{lem:CLT md}
For each $n\geq1$, let $\left\{Y_{ni}, \cF_{ni}\right\}_{i\geq1}$ be a martingale difference array on $(\Omega,\cF, \P)$, with $\E|Y_{ni}|^2<\infty$ for all $n\geq1$ and let $\tau_n$ be a finite stopping time w.r.t. $\cF_{ni}$. Suppose that for some constant $\sigma^2>0$, as $n\rightarrow\infty$,
\begin{eqnarray*}
\sum_{i=1}^{\tau_n} \E\left[\left. |Y_{ni}|^2\right|\cF_{ni}\right]\to\sigma^2, \ in\ probability,
\end{eqnarray*}
and that for each $\varepsilon>0$,
\begin{eqnarray*}
\sum_{i=1}^{\tau_n} \E\left[\left. |Y_{ni}|^2\mathds{1}_{\left\{\left|Y_{ni}\right|>\varepsilon\right\}}\right|\cF_{ni}\right]\to0, \ in\ probability,
\end{eqnarray*}
Then,
\begin{eqnarray*}
\sum_{i=1}^{\tau_n} Y_{ni}\Rightarrow N(0,\sigma^2), \ as\ n\rightarrow\infty.
\end{eqnarray*}
\end{lemma}

\noindent {\bf Proof of Theorem \ref{thm:CLT tilde M}.} By (\ref{eqn:tilde M decomp}),
\begin{eqnarray}\label{tilde M CLT0}
\sqrt{n}\left(\widetilde M_n-\mu^*\right)=\frac{Z_n}{\sqrt{n}}+\sqrt{n}\left(\frac{1}{n}\sum_{j=1}^n \mu_{I_j}-\mu^*\right).
\end{eqnarray}

We analyze the two terms on the RHS of (\ref{tilde M CLT0}) separately. We have shown that $Z_n$ is a martingale. So $Z_n-Z_{n-1}=X_{I_n,n}-\mu_{I_n}$ is a martingale difference array (see Definition \ref{def:md}). To show the asymptotic normality of $Z_n/\sqrt{n}$, it suffices to prove that the two conditions of Lemma \ref{lem:CLT md} are satisfied for $Y_{nj}=\left(X_{I_j,j}-\mu_{I_j}\right)/\sqrt{n}$, $\tau_n=n$ and $\cF_{nj}=\cF_j$.

First, note that
\begin{eqnarray*}
&&\sum_{j=1}^n \E\left[\left. \left|\left(X_{I_j,j}-\mu_{I_j}\right)/\sqrt{n}\right|^2\right|\cF_{j-1}\right]=\frac{1}{n}\sum_{j=1}^n \E\left[\left. X_{I_j,j}^2-2\mu_{I_j}X_{I_j,j}+\mu_{I_j}^2\right|\cF_{j-1}\right]\\
&=&\frac{1}{n}\sum_{j=1}^n \left\{\E\left[\left. X_{I_j,j}^2 \right|\cF_{j-1}\right]-2\mu_{I_j}\E\left[\left.X_{I_j,j}\right|\cF_{j-1}\right]+\mu_{I_j}^2\right\}=\frac{1}{n}\sum_{j=1}^n \left\{\E\left[\left. X_{I_j,j}^2 \right|\cF_{j-1}\right]-\mu_{I_j}^2\right\}\\
&=&\frac{1}{n}\sum_{j=1}^n \sigma_{I_j}^2=\frac{1}{n}\sum_{k=1}^K T_k(n) \sigma_k^2 \overset{a.s.}{\longrightarrow} \sigma_{k^*}^2,\ {\rm as}\ n\to\infty,
\end{eqnarray*}
where the convergence follows from Theorem \ref{thm:Tk a.s.}.

Second, by Cauchy-Schwarz inequality,
\begin{align*}
&\sum_{j=1}^n \E\left[\left. \left|\left(X_{I_j,j}-\mu_{I_j}\right)/\sqrt{n}\right|^2\1\left\{\left|\left(X_{I_j,j}-\mu_{I_j}\right)/\sqrt{n}\right|>\varepsilon\right\}\right|\cF_{j-1}\right]\\
\leq&\sum_{j=1}^n \left(\E\left[\left. \left|\left(X_{I_j,j}-\mu_{I_j}\right)/\sqrt{n}\right|^4\right|\cF_{j-1}\right]\right)^\frac{1}{2}
\left(\E\left[\left. \1\left\{\left|\left(X_{I_j,j}-\mu_{I_j}\right)/\sqrt{n}\right|>\varepsilon\right\}\right|\cF_{j-1}\right]\right)^\frac{1}{2}\\
=&\frac{1}{n}\sum_{j=1}^n \left(\E\left[\left. \left|X_{I_j,j}-\mu_{I_j}\right|^4\right|\cF_{j-1}\right]\right)^\frac{1}{2}
\P\left(\left. \left|\left(X_{I_j,j}-\mu_{I_j}\right)\right|>\varepsilon\sqrt{n}\right|\cF_{j-1}\right)^\frac{1}{2}\\
=&\frac{1}{n}\sum_{j=1}^n \vartheta_{I_j}^2
\P\left(\left. \left|\left(X_{I_j,j}-\mu_{I_j}\right)\right|>\varepsilon\sqrt{n}\right|\cF_{j-1}\right)^\frac{1}{2}
\leq\frac{1}{n}\sum_{j=1}^n \vartheta_{I_j}^2
\left((\varepsilon^2 n)^{-1}\E\left[\left. \left|X_{I_j,j}-\mu_{I_j}\right|^2\right|\cF_{j-1}\right]\right)^\frac{1}{2}\\
=& \frac{1}{n}\sum_{j=1}^n \vartheta_{I_j}^2
\left((\varepsilon^2 n)^{-1}\sigma_{I_j}^2\right)^\frac{1}{2}
\leq \frac{\bar\vartheta^2 \bar\sigma}{\varepsilon\sqrt{n}} \to 0, \mbox{ as } n\to\infty,
\end{align*}
where $\vartheta_k^4=\E[|X_k-\mu_k|^4]$, $\bar\vartheta^4 = \max\{\vartheta_1^4,...,\vartheta_K^4\}$ and $\bar\sigma^2 = \max\{\sigma_1^2,...,\sigma_K^2\}$.

Therefore, by Lemma \ref{lem:CLT md}, we have
\begin{eqnarray*}
Z_n/\sqrt{n}\Rightarrow N(0,\sigma_{k^*}^2), \ \ as\ \ n\rightarrow\infty.
\end{eqnarray*}

It remains to prove that the second term on the RHS of (\ref{tilde M CLT0}) converges to 0 in probability. Note that
\begin{eqnarray*}
\frac{1}{n}\sum_{j=1}^n \mu_{I_j}-\mu_{k^*}&=&\frac{1}{n}\sum_{k=1}^K T_k(n) \mu_k -\mu_{k^*}=\sum_{k\neq k^*}\frac{T_k(n)}{n} \mu_k + \left(\frac{T_{k^*}(n)}{n}-1\right)\mu_{k^*}\\
&=&\sum_{k\neq k^*}\frac{T_k(n)}{n} \mu_k - \sum_{k\neq k^*}\frac{T_k(n)}{n}\mu_{k^*}=\sum_{k\neq k^*}\frac{T_k(n)}{n} \left(\mu_k - \mu_{k^*}\right).
\end{eqnarray*}
Therefore, the first moment of the second term on the RHS of (\ref{tilde M CLT0}) satisfies
\begin{eqnarray*}
\E\left| \sqrt{n}\left(\frac{1}{n}\sum_{j=1}^n \mu_{I_j}-\mu_{k^*}\right) \right|
=\E\left| \sum_{k\neq k^*}\frac{T_k(n)}{\sqrt{n}} \left(\mu_k - \mu_{k^*}\right)  \right|
\leq \sum_{k\neq k^*}\frac{\E T_k(n)}{\sqrt{n}} \Delta_k \to 0.
\end{eqnarray*}
where the convergence follows from Proposition \ref{prop:mom of Tk} with $p=1$ and the condition $\nu_n\in\left[\alpha\log n, n^{\frac{1}{2}-\delta}\right]$  $(0<\delta<1, \alpha\geq4\bar\gamma^2)$.

Applying Slutsky's Theorem to the RHS of (\ref{tilde M CLT0}) leads to the conclusion.

\subsection{Proof of Proposition \ref{prop:CI tilde M}}\label{sec:ProofSigmaConsistency}

Note that
\begin{eqnarray}\label{sigma k}
\sigma_k^2=\E X_k^2 - \left(\E X_k\right)^2.
\end{eqnarray}
We first prove that $\hat\sigma_k^2$ converges to $\sigma_k^2$ in probability. From Theorem 2 of Lai and Robbins (1985), it follows that for $k\neq k^*$, $T_k(n)\to\infty$ in probability; and from Theorem \ref{thm:Tk a.s.}, it follows that $T_{k^*}(n)\to\infty$ a.s. By Theorem 2 in Richter (1965), we have that for $k=1,...,K$,
\begin{eqnarray*}
\frac{1}{T_k(n)}\sum_{j=1}^{T_k(n)}X_{k,j}^2\overset{\P}{\longrightarrow}\E X_k^2,\ \ {\rm and}\ \ \frac{1}{T_k(n)}\sum_{j=1}^{T_k(n)}X_{k,j}\overset{\P}{\longrightarrow}\E X_k.
\end{eqnarray*}
By the continuous mapping theorem, it follows that $\hat\sigma_k^2\to\sigma_k^2$ in probability.

Then by Theorem \ref{thm:Tk a.s.} and the continuous mapping theorem, $\widetilde\sigma_n^2$ converges to $\sigma_{k^*}^2$ in probability.

\subsection{Proof of Theorem \ref{TIstar a.s.}}\label{sec:ProofConsistencyLSA}

By Theorem \ref{thm:Tk a.s.}, we have
\begin{eqnarray}\label{Tk a.s.infty}
T_{k^*}(n)\overset{a.s.}{\longrightarrow}\infty.
\end{eqnarray}

By the strong law of large numbers (SLLN) and Theorem 1 in Richter (1965) or Theorem 8.2 in Chapter 6 of Gut (2013), we have
\begin{eqnarray*}
\frac{1}{T_{k^*}(n)}\sum_{j=1}^{T_{k^*}(n)}X_{k^*,j}\overset{a.s.}{\longrightarrow} \mu^*,\ \ {\rm as}\ \ n\to\infty.
\end{eqnarray*}
Denote $M_k=\frac{1}{T_{k}(n)}\sum_{j=1}^{T_{k}(n)}X_{k,j}$.
The definition of almost sure convergence says that
\begin{eqnarray}\label{Prob a.s. Tk star}
\P\left(\omega\in\Omega: \forall 0<\varepsilon<1, \exists N'=N'(\omega,\varepsilon)\in \N, s.t.\ \forall n>N', \left|M_{k^*}-\mu^*\right|<\varepsilon\right)=1.
\end{eqnarray}
From (\ref{a.s. In0}), it follows that
\begin{eqnarray}\label{Prob a.s. TIn}
\P\left(\omega\in\Omega: \forall 0<\varepsilon<1, \exists N=N(\omega,\varepsilon)\in \N, s.t.\ \forall n>N, M_{I^*_n}=M_{k^*}\right)=1.
\end{eqnarray}

Combining (\ref{Prob a.s. Tk star}) and (\ref{Prob a.s. TIn}) leads to
\begin{eqnarray*}
\P\left(\omega\in\Omega: \forall 0<\varepsilon<1, \exists \bar N=\max\{N,N'\}\in \N, s.t.\ \forall n>\bar N, \left|M_{I^*_n}-\mu^*\right|<\varepsilon\right)=1,
\end{eqnarray*}
implying that
\begin{eqnarray*}
M_{I^*_n}\overset{a.s.}{\rightarrow} \mu^*,\ \ {\rm as}\ \ n\to\infty.
\end{eqnarray*}

\subsection{Proof of Theorem \ref{thm:MSE Mstar}}\label{sec:ProofMSE_LSA}

For notational ease, denote
\[Z_{n,k} = \sum_{j=1}^{T_k(n)} \left( X_{k,j} - \mu_k \right),\quad k = 1,\ldots,K.\]

Recall that sub-Gaussian distributions have finite moments up to any order. By Assumption \ref{ass:subGaussian} and the fact that $T_k(n)\le n$, it can be seen that $\{X_{k,j}, j\ge 1\}$ satisfies the conditions in Lemma \ref{lem:SumBound} for any $k=1,\ldots,K$. We first note the following two facts that are useful in the proof. For any positive integer $p$, by Proposition \ref{prop:mom of Tk},
\begin{eqnarray}\label{eqn:Fact1}
\Pr(I_n^*\neq k^*) \le \sum_{k\neq k^*}\Pr(I_n^*=k) \le \sum_{k\neq k^*} \Pr(T_k(n)\ge n/K) \le K^{p+1} \E\left[ T_k(n)^p \right]/n^p = O\left({\nu_n^p\over n^p} + {1\over n^{3}} \right),
\end{eqnarray}
and by Lemma \ref{lem:SumBound},
\begin{eqnarray}\label{eqn:Fact2}
\E\left[ \left|Z_{n,k}\right|^p \right] = O(n^{p/2}),\quad k=1,\ldots, K. 
\end{eqnarray}

It suffices to analyze the bias and variance of $M_{I_n^*}$, and then the MSE follows straightforwardly. Consider the bias first. Note that
\begin{eqnarray}
\lefteqn{{\rm Bias}(M_{I_n^*}) = \E\left[M_{I_n^*} - \mu_{k^*}\right] = \E\left[\frac{1}{T_{I_n^*}(n)} \sum_{j=1}^{T_{I_n^*}(n)} (X_{I_n^*,j} - \mu_{k^*})\right]}\nonumber\\
&=& \E\left[\frac{Z_{n,k^*}}{T_{k^*}(n)}\1\{I_n^*=k^*\}\right] + \sum_{k\neq k^*} \E\left[\frac{Z_{n,k^*}}{T_k(n)} \1\{I_n^*=k\}\right]\nonumber\\
&=&\E\left[\frac{Z_{n,k^*}}{T_{k^*}(n)}\right] - \sum_{k\neq k^*} \E\left[\frac{Z_{n,k^*}}{T_{k^*}(n)} \1\{I_n^*=k\}\right] + \sum_{k\neq k^*} \E\left[\frac{Z_{n,k^*}}{T_k(n)} \1\{I_n^*=k\}\right]\label{eqn:TempNum}\\
&=& Q_1 - Q_2 + Q_3,\nonumber
\end{eqnarray}
where $Q_1, Q_2$ and $Q_3$ denote the three terms on the RHS of (\ref{eqn:TempNum}), respectively. 

Recall that $\E[Z_{n,k^*}]=0$ due to Wald's equation. It follows that
\begin{align}
Q_1 =& \E\left[\frac{Z_{n,k^*}}{T_{k^*}(n)} - \frac{Z_{n,k^*}}{\E T_{k^*}(n)}\right] = \E\left[\frac{Z_{n,k^*}}{T_{k^*}(n)}\left(1 - \frac{T_{k^*}(n)}{\E T_{k^*}(n)}\right)\right]
= \Cov\left[\frac{Z_{n,k^*}}{T_{k^*}(n)}, 1 - \frac{T_{k^*}(n)}{\E T_{k^*}(n)}\right]\nonumber\\
=& -\Cov\left[\frac{Z_{n,k^*}}{T_{k^*}(n)}, \frac{T_{k^*}(n)}{\E T_{k^*}(n)}\right]
= -\rho \cdot \frac{\sqrt{\Var[Z_{n,k^*}/T_{k^*}(n)] \Var[T_{k^*}(n)]}}{\E T_{k^*}(n)},\label{LSAbias Q1}
\end{align}
where $\rho$ denote the correlation coefficient between $Z_{n,k^*}/T_{k^*}(n)$ and $T_{k^*}(n)$, i.e.
\begin{align*}
\rho=\frac{ \Cov\left[Z_{n,k^*}/T_{k^*}(n),T_{k^*}(n)\right]}{\sqrt{\Var[Z_{n,k^*}/T_{k^*}(n)] \Var[T_{k^*}(n)]}}.
\end{align*}
%Then,
%\begin{align}
%Q_1 =  -\rho \cdot \frac{\sqrt{\Var[Z_{n,k^*}/T_{k^*}(n)] \Var[T_{k^*}(n)]}}{\E T_{k^*}(n)}.
%\end{align}
It should also be noted that $\rho>0$ (and thus $Q_1<0$), because $Z_{n,k^*}/T_{k^*}(n)$ and $T_{k^*}(n)$ are positive correlated. To see this, note that
\begin{align*}
&\Cov\left[\frac{Z_{n,k^*}}{T_{k^*}(n)},T_{k^*}(n)\right] = \Cov\left[\frac{\bar X_{k^*}[T_{k^*}(n)]T_{k^*}(n) }{T_{k^*}(n)}-\mu_{k^*},T_{k^*}(n)\right] = \Cov\left[\bar X_{k^*}[T_{k^*}(n)],T_{k^*}(n)\right]>0,
\end{align*}
where the inequality follows naturally from the way that GUCB is defined, i.e., larger sample mean leads to higher chance of sampling from arm $k^*$. 

By (\ref{LSAbias Q1}), to understand the order $Q_1$, we only need to analyze $\E T_{k^*}(n)$, $\Var[T_{k^*}(n)]$ and $\Var[Z_{n,k^*}/T_{k^*}(n)]$. Note that by Proposition \ref{prop:mom of Tk},
\begin{align}\label{LSAbias Q101}
n > \E T_{k^*}(n) = n - \sum_{k\neq k^*} \E T_k(n) \geq n - \sum_{k\neq k^*}\left( \frac{8\nu_n}{\Delta_k^2} + 5\right) \geq n - Kc_1\nu_n,
\end{align}
and
\begin{align}\label{LSAbias Q102}
\Var[T_{k^*}(n)] =& \Var\left[n-\sum_{k\neq k^*}T_k(n)\right] = \Var\left[\sum_{k\neq k^*}T_k(n)\right]
\leq \E\left[\left(\sum_{k\neq k^*}T_k(n)\right)^2\right]\nonumber \\
\leq& (K-1)\E\left[\sum_{k\neq k^*}T_k(n)^2\right]
\leq K\sum_{k\neq k^*} \left( \frac{8\nu_n}{\Delta_k^2} + 5 \right)^2
\leq K^2 c_2^2 \nu_n^2,
\end{align}
where the third inequality follows from Proposition \ref{prop:mom of Tk} with $p=2$.

Moreover,
\begin{align}\label{LSAbias Q103}
\Var\left[\frac{Z_{n,k^*}}{T_{k^*}(n)}\right] \leq \E\left[\frac{Z_{n,k^*}^2}{T_{k^*}^2(n)}\right] = \E\left[\frac{Z_{n,k^*}^2}{T_{I_n^*}^2(n)}\1\{I_n^* = k^*\}\right] + \E\left[\frac{Z_{n,k^*}^2}{T_{k^*}(n)^2}\1\{I_n^* \neq k^*\}\right],
\end{align}
We analyze the RHS of (\ref{LSAbias Q103}). Because $T_{I_n^*}\ge n/K$, 
\begin{eqnarray}\label{eqn:TempR3}
\E\left[\frac{Z_{n,k^*}^2}{T_{I_n^*}^2(n)}\1\{I_n^* = k^*\}\right]\le {K^2\over n^2}\E\left[ Z_{n,k^*}^2 \right] = {\sigma_{k^*}^2 K^2\over n^2}\E\left[ T_{k^*}(n) \right]\le {\sigma_{k^*}^2 K^2\over n},
\end{eqnarray}
where the first equality follows from Wald's lemma (Theorem 13.2.14 in Athreya and Lahiri 2006).

We also note that
\begin{eqnarray}
\lefteqn{\E\left[\frac{Z_{n,k^*}^2}{T_{k^*}(n)^2}\1\{I_n^* \neq k^*\}\right]\le \E\left[Z_{n,k^*}^2\1\{I_n^* \neq k^*\}\right]}\nonumber\\
&\le & \E^{p_1-1\over p_1}\left( \1\{I_n^* \neq k^*\} \right)\E^{1\over p_1}\left( Z_{n,k^*}^{2p_1} \right)\nonumber\\
& = & O\left( \left({\nu_n^{p_2}\over n^{p_2}} + {1\over n^3} \right)^{p_1-1\over p_1}\cdot n^{{2p_1\over 2}\cdot {1\over p_1}} \right) = O\left( \left({n^{(1-\delta)p_2(p_1-1)/p_1}\over n^{p_2(p_1-1)/p_1}} + {1\over n^{3(p_1-1)/p_1}} \right)\cdot n \right) = o\left({1\over n} \right),\label{eqn:TempR2}
\end{eqnarray}
where the second inequality follows from H\"older's inequality, and the third-to-last equality is due to the facts in (\ref{eqn:Fact1}) and (\ref{eqn:Fact2}) with $p_1$ and $p_2$ satisfying $p_1=3/(1-\epsilon)$ and $p_2>3/\delta$ for some $\epsilon\in (0,1)$. 

Combining Equations (\ref{LSAbias Q1})-(\ref{eqn:TempR2}) yields 
\begin{eqnarray*}
|Q_1|\le {K^2\sigma_{k^*}^2 c_2\nu_n + o(1)\over \sqrt{n}\left( n - Kc_1\nu_n \right)} = O\left( {\nu_n\over n^{3/2}} \right).
\end{eqnarray*}

Using the same argument in proving (\ref{eqn:TempR2}), we have
\begin{eqnarray*}
\lefteqn{|Q_2| \le \sum_{k\neq k^*} \E\left[\left|Z_{n,k^*}\right| \1\{I_n^*=k\}\right]}\\
&\le & (K-1) \E^{p_1-1\over p_1}\left( \1\{I_n^* \neq k^*\} \right)\E^{1\over p_1}\left| Z_{n,k^*}^{p_1} \right|\nonumber\\
&=& O\left( \left({n^{(1-\delta)p_2(p_1-1)/p_1}\over n^{p_2(p_1-1)/p_1}} + {1\over n^{3(p_1-1)/p_1}} \right)\cdot n^{1\over 2} \right)=o\left( {\nu_n\over n^{3/2}} \right)
\end{eqnarray*}
with $p_1$ and $p_2$ satisfying $p_1=3/(1-\epsilon)$ and $p_2>3/(2\delta)$ for some $\epsilon\in (0,1)$. 

In a similar manner, it can be shown that $|Q_3|=o\left( {\nu_n\over n^{3/2}} \right)$, and the details are omitted. Therefore, by (\ref{eqn:TempNum}), 
\[\left|{\rm Bias}\left( M_{I_n^*} \right)\right| \le {K^2\sigma_{k^*}^2 c_2\nu_n + c_3\over \sqrt{n}\left( n - Kc_1\nu_n \right)} +o\left( {\nu_n\over n^{3/2}} \right) = O\left( {\nu_n\over n^{3/2}} \right). \]

Note further that the order of $|Q_1|$ dominates those of $|Q_2|$ and $|Q_3|$. Because $Q_1<0$, we know that ${\rm Bias}\left( M_{I_n^*} \right)<0$ for sufficiently large $n$. 

We now consider the variance. Note that
\begin{align}\label{LSAVar0}
\Var\left[M_{I_n^*}\right] \leq& \E\left[\left|M_{I_n^*}-\mu_{k^*}\right|^2\right]\nonumber \\
=& \E\left[\left|M_{I_n^*}-\mu_{k^*}\right|^2\1\{I_n^*=k^*\}\right] + \sum_{k\neq k^*}\E\left[\left|M_{I_n^*}-\mu_{k^*}\right|^2\1\{I_n^*=k\}\right]\nonumber \\
=& \E\left[\frac{Z_{n,k^*}^2}{T_{k^*}^2(n)}\1\{I_n^*=k^*\}\right] + \sum_{k\neq k^*}\E\left[{ (Z_{n,k} - \Delta_k)^2 \over T_k^2(n)}\1\{I_n^*=k\}\right]\nonumber \\
\leq& \E\left[\frac{Z_{n,k^*}^2}{T_{k^*}^2(n)}\1\{I_n^*=k^*\}\right] + 2\sum_{k\neq k^*}\E\left[{ Z_{n,k}^2 \over T_k^2(n)}\1\{I_n^*=k\}\right]+ 2\sum_{k\neq k^*}\Delta_k^2\E\left[{1 \over T_k^2(n)}\1\{I_n^*=k\}\right].
\end{align}
We note that by (\ref{eqn:TempR3}), 
\[\E\left[\frac{Z_{n,k^*}^2}{T_{k^*}^2(n)}\1\{I_n^* = k^*\}\right]\le  {\sigma_{k^*}^2 K^2\over n}.\]Similar to the proof of (\ref{eqn:TempR2}), we have
\[ 2\sum_{k\neq k^*}\E\left[{ Z_{n,k}^2 \over T_k^2(n)}\1\{I_n^*=k\}\right] = o\left({1\over n} \right).\]Moreover, using the fact in (\ref{eqn:Fact1}),
\[2\sum_{k\neq k^*}\Delta_k^2\E\left[{1 \over T_k^2(n)}\1\{I_n^*=k\}\right] \le 2\sum_{k\neq k^*}\Delta_k^2\E\left[\1\{I_n^*=k\}\right]  = O\left( {\nu_n^p\over n^p} + {1\over n^3} \right)=O\left( {n^{(1-\delta)p}\over n^p} + {1\over n^3} \right) = o\left( {1\over n}\right),\]with $p>1/\delta$. 

Therefore, 
\[\Var\left[M_{I_n^*}\right]  \le {\sigma_{k^*}^2 K^2\over n} + o\left( {1\over n}\right),\]which completes the proof.

\subsection{Proof of Theorem \ref{thm:CLT Mstar}}\label{sec:ProofCLT_LSA}

Let $M_k=\frac{1}{T_{k}(n)}\sum_{j=1}^{T_{k}(n)}X_{kj}$. By (\ref{Prob a.s. TIn}),
\begin{eqnarray*}
\P\left(\omega\in\Omega: \forall 0<\varepsilon<1, \exists N=N(\omega,\varepsilon)\in \N, s.t.\ \forall n>N, \sqrt{n}M_{I^*_n}=\sqrt{n}M_{k^*}\right)=1,
\end{eqnarray*}
i.e.,
\begin{eqnarray*}
\sqrt{n}M_{I^*_n}-\sqrt{n}M_{k^*}\overset{a.s.}{\longrightarrow} 0,\ \ {\rm as}\ \ n\to\infty.
\end{eqnarray*}
Note that
\begin{eqnarray*}
\sqrt{n}\cdot\left(M_{I^*_n}-\mu^*\right)=\sqrt{n}M_{I^*_n}-\sqrt{n}M_{k^*}+\sqrt{n}\cdot\left(M_{k^*}-\mu^*\right).
\end{eqnarray*}
By Slutsky Theorem, to prove the theorem, it suffices to show
\begin{eqnarray}\label{CLT sum Tk}
\sqrt{n}\cdot\left(M_{k^*}-\mu^*\right) \Rightarrow N(0,\sigma_{k^*}^2),\ \ {\rm as}\ \ n\to\infty,
\end{eqnarray}
which follows immediately from (\ref{Tk a.s.infty}) and Theorem 1 of R\'{e}nyi (1957) or Theorem 3.2 in Chapter 7 of Gut (2013), and the proof is completed.

\subsection{Proof of Proposition \ref{prop:CI Mstar}}\label{sec:ProofConsistencyVar_LSA}

Note that $\sigma_{k^*}^2=\E X_{k^*}^2 - \E^2 X_{k^*}$. It suffices to prove the two terms on the RHS of (\ref{bar sigma}) converge to two terms of $\sigma_{k^*}^2$, respectively.

By Theorem \ref{TIstar a.s.} and the continuous mapping theorem, the second term on the RHS of (\ref{bar sigma}) converges almost surely to $\E^2 X_{k^*}$. Similar to the proof of Theorem \ref{TIstar a.s.}, we can show that the first term on the RHS of (\ref{bar sigma}) converge to $\E X_{k^*}^2$ almost surely, thus completing the proof.

\section*{References}

\begin{hangref}
\item Agrawal, R. 1995. Sample mean based index policies by $o(\log n)$ regret for the multi-armed bandit problem. {\it Advances in Applied Probability}, {\bf 27(4)}, 1054-1078.

%\item Artzner, P., F. Delbaen, J. M. Eber, D. Heath. 1999. Coherent measures of risk. {\it Mathematical Finance}, {\bf 9(3)}, 203-228.

\item Athreya, K. B., S. N. Lahiri. 2006. {\it Measure Theory and Probability Theory}. Springer Science \& Business Media.

\item Auer, P.,  N. Cesa-Bianchi, P. Fischer. 2002. Finite-time analysis of the multiarmed bandit problem. {\it Machine Learning}, {\bf 47(2-3)}, 235-256.

%\item Blanchet, J., K. Murthy. 2019. Quantifying distributional model risk via optimal transport. {\it Mathematics of Operations Research}, {\bf 44(2)}, 565-600.

\item Chang, H. S., M. C. Fu, J. Hu, S. I. Marcus. 2005. An adaptive sampling algorithm for solving Markov decision processes. {\it Operations Research}, {\bf 53(1)}, 126-139.

\item Chen, H. J., E. J. Dudewicz. 1976. Procedures for fixed-width interval estimation of the largest normal mean. {\it Journal of the American Statistical Association}, {\bf 71(355)}, 752-756. 

%\item Chung, K. L. 2001. {\it A Course on Probability Theory}, Academic Press, New York.

%\item Coulom, R. 2006. Efficient selectivity and backup operators in Monte-Carlo tree search. In {\it International Conference on Computers and Games}, 72-83. Springer, Berlin, Heidelberg.

%\item Drzisga, D., B. Gmeiner, U. R\"{u}de, R. Scheichl, B. Wohlmuth. 2017. Scheduling massively parallel multigrid for multilevel Monte Carlo methods. {\it SIAM Journal on Scientific Computing}, {\bf 39(5)}, 873-897.

%\item Dean, B. C., M. X. Goemans, J. Vondr\'ak. 2008. Approximating the stochastic knapsack problem: The benefit of adaptivity. {\it Mathematics of Operations Research}, {\bf 33(4)}, 945-964.

\item D'Eramo, C., A. Nuara, M. Restelli. 2016. Estimating the maximum expected value through Gaussian approximation. {\it International Conference on Machine Learning}, 1032-1040.

\item Dmitrienko, A., J. C. Hsu. 2006. Multiple testing in clinical trials. S. Kotz, C. B. Read, N. Balakrishnan, B. Z. Vidakovic (Eds.), Encyclopedia of Statistical Sciences. Wiley, New Jersey. %5111-5117

%\item Durrett, R. 2019. {\it Probability: Theory and Examples}, (Vol.49). Cambridge University Press.

%\item Edelkamp, S., M. Gath, C. Greulich, M. Humann, O. Herzog, M. Lawo. 2016. Monte-Carlo tree search for logistics. In {\it Commercial Transport}, 427-440. Springer, Cham.

%\item Dzyabura D., J. R. Hauser. 2019. Recommending products when consumers learn their preference weights. {\it Marketing Science}, {\bf 38(3)}, 417-441.

\item Fan, W., L. J. Hong, B. L. Nelson. 2016. Indifference-zone-free selection of the best. {\it Operations Research}, {\bf 64(6)}, 1499-1514.

\item FDA (Food and Drug Administration). 2019. Adaptive designs for clinical trials of drugs and biologics: Guidance for industry. Available at https://www.fda.gov/media/78495/download.

\item Fu, M. C. 2017. Markov decision processes, AlphaGo, and Monte Carlo tree search: Back to the future. {\it INFORMS TutORials in Operations Research}, 68-88.

%\item Ghosh, S., H. Lam. 2019. Robust analysis in stochastic simulation: Computation and performance guarantees. {\it Operations Research}, {\bf 67(1)}, 232-249.

\item Gut, A. 2013. {\it Probability: A Graduate Course}. Springer Science \& Business Media.

%\item Helm, J. E., M. S. Lavieri, M. P. Van Oyen, J. D. Stein, D. C. Musch. 2015. Dynamic forecasting and control algorithms of glaucoma progression for clinician decision support. {\it Operations Research}, {\bf 63(5)}, 979-999.

\item Hao, B., Y. Abbasi-Yadkori, Z. Wen, G. Cheng. 2019. Bootstrapping upper confidence bound. {\it Advances in Neural Information Processing Systems}, 12123-12133. 

%\item Hong, L. J., G. Liu. 2009. Simulating sensitivities of conditional value at risk. {\it Management Science}, {\bf 55(2)}, 281-293.

%\item Kaufmann, E., W. M. Koolen, A. Garivier. 2018. Sequential test for the lowest mean: From Thompson to Murphy sampling. {\it Adcances in Neural Information Processing Systems (NIPS)}, 6333-6343.

\item Kim, S.-H., B. L. Nelson. 2006. Selecting the best system. In {\it Handbook in OR \& MS: Simulation}, S. G. Henderson and B. L. Nelson (Eds.), 501-534, Elsevier Science.

%\item Kleywegt A. J., J. D. Papastavrou. 2001. The dynamic and stochastic knapsack problem with random sized items. {\it Operations research}, {\bf 49(1)}, 26-41.

\item Kocsis, L., C. Szepesv\'ari. 2006. Bandit based Monte-Carlo planning. {\it Machine Learning: ECML}, 282-293.

\item Lai, T. L., H. Robbins. 1985. Asymptotically efficient adaptive allocation rules. {\it Advances in Applied Mathematics}, {\bf 6(1)}, 4-22.

\item Lesnevski, V., B. L. Nelson, J. Staum. 2007. Simulation of coherent risk measures based on generalized scenarios. {\it Management Science}, {\bf 53(11)}, 1756-1769.

%\item Mandelbaum, A. 2014. Service Engineering (096324). Winter 2014. http://ie.technion.ac.il/serveng/, accessed 2020-07-13.

\item Mozgunov, P., T. Jaki. 2020. An information theoretic approach for selecting arms in clinical trials. {\it Journal of the Royal Statistical Society Series B}, {\bf 82(5)}, 1223-1247.

\item Qin, C., D. Klabjan, D. Russo. 2017. Improving the expected improvement algorithm. {\it Proceedings of the 31st International Conference on Neural Information Processing Systems}, 5387-5397.

\item R\'{e}nyi, A. 1957. On the asymptotic distribution of the sum of a random number of independent random variables. {\it Acta Mathematica Academiae Scientiarum Hungarica}, {\bf 8(1-2)}, 193-199.

\item Richter, W. 1965. Limit theorems for sequences of random variables with sequences of random indeces. {\it Theory of Probability \& Its Applications}, {\bf 10(1)}, 74-84.

%\item Rimmel, A., F. Teytaud, T. Cazenave. 2011. Optimization of the nested Monte-Carlo algorithm on the traveling salesman problem with time windows. In {\it European Conference on the Applications of Evolutionary Computation}, 501-510. Springer, Berlin, Heidelberg.

\item Robertson, D. S., J. M. S. Wason. 2019. Familywise error control in multi-armed response-adaptive trials. {\it Biometrics}, {\bf 75}, 885-894. 

%\item Ross, S. M. 1996. {\it Stochastic Processes}, 2nd Edition. John Wiley, New York.

\item Russo, D. 2020. Simple bayesian algorithms for best arm identification. {\it Operations Research}, {\bf 68(6)}, 1625-1647.
    
\item Shamir, O. 2011. A variant of Azuma's inequality for martingales with subgaussian tails. arXiv preprint arXiv:1110.2392.

%\item Searle, S. R., G. Casella, C. E. McCulloch. 2009. {\it Variance Components}, (Vol.391). John Wiley \& Sons.

%\item Shi, W., X. Chen. 2018. Efficient budget allocation strategies for elementary effects method in stochastic simulation. {\it Naval Research Logistics}, {\bf 65(3)}, 218-241.

%\item Shin, D., M. Broadie, A. Zeevi. 2018. Tractable sampling strategies for ordinal optimization. {\it Operations Research}, {\bf 66(6)}, 1693-1712.

%\item Shapiro A., D. Dentcheva, A. Ruszczynski. 2014. {\it Lectures on Stochastic Programming: Modeling and Theory}, (Vol. 16). SIAM, Philadelphia.

%\item Silver, D., A. Huang, C. J. Maddison, A. Guez, L. Sifre, G. Van Den Driessche, J. Schrittwieser, I. Antonoglou, V. Panneershelvam, M. Lanctot, S. Dieleman, D. Grewe, J. Nham, N. Kalchbrenner, I. Sutskever, T. Lillicrap, M. Leach, K. Kavukcuoglu, T. Graepel, D. Hassabis. 2016. Mastering the game of Go with deep neural networks and tree search. {\it Nature}, {\bf 529(7587)}, 484-489.

%\item Silver, D., T. Hubert, J. Schrittwieser, I. Antonoglou, M. Lai, A. Guez, M. Lanctot, L. Sifre, D. Kumaran, T. Graepel, T. Lillicrap, K. Simonyan, D. Hassabis. 2018. A general reinforcement learning algorithm that masters chess, shogi, and Go through self-play. {\it Science}, {\bf 362(6419)}, 1140-1144.

\item Smith, J. E., R. L. Windler. 2006. The optimizer's curse: Skepticism and postdecision surprise in decision analysis. {\it Management Science}, {\bf 52(3)}, 311-322.

%\item Sun, Y., D. W. Apley, J. Staum. 2011. Efficient nested simulation for estimating the variance of a conditional expectation. {\it Operations Research}, {\bf 59(4)}, 998-1007.

%\item Vercraene, S., J. P. Gayon, F. Karaesmen. 2017. Effects of System Parameters on the Optimal Cost and Policy in a Class of Multidimensional Queueing Control Problems. {\it Operations Research}, {\bf 66(1)}, 150-162.

\item van Hasselt, H. 2010. Double $q$-learning. {\it Advances in Neural Information Processing Systems (NIPS)}, 2613-2621.

\item van Hasselt, H. 2013. Estimating the maximum expected value: an analysis of (nested) cross validation
and the maximum sample average. Available at arXiv preprint arXiv:1302.7175.

\item van Hasselt, H.,, A. Guez, D. Silver. 2016. Deep reinforcement learning with double Q-learning.
{\it Proceedings of the AAAI Conference on Artificial Intelligence}, {\bf 30}, 2094-2100.

\item Villar, S. S., J. Bowden, J. Wason. 2015. Multi-armed bandit models for the optimal design of clinical trials: Benefits and challenges. {\it Statistical Science}, {\bf 30(2)}, 199.

\item Wainwright, M. J. 2019. High-dimensional Statistics: A Non-asymptotic Viewpoint (Vol. 48). Cambridge University Press.

%\item Xu, M., T. Qin, T.-Y. Liu. 2013. Estimation bias in multi-armed bandit algorithms for search advertising. {\it Advances in Neural Information Processing Systems (NIPS)}, 2400-2408.

%\item Zhong, Y., Hong, L. J., Liu, G. 2021. Earning and learning with varying cost. {\it Production and Operations Management}, {\bf 30(8)}, 2379-2394.
\end{hangref}

\end{document}